\documentclass[12pt,reqno]{article}
\usepackage{amssymb}
\usepackage{amsmath}
\usepackage{amsthm}
\usepackage{amsfonts}
\usepackage{url}
\usepackage{hyperref}
\usepackage{breakurl}
\usepackage{rotating}
\usepackage{adjustbox}
\def\modd#1 #2{#1\ \mbox{\rm (mod}\ #2\mbox{\rm )}}
\usepackage{pdflscape}
\usepackage{microtype}
\newcommand{\comment}[1]{}
\usepackage{caption}
\newcommand{\BE}{\begin{equation}}
\newcommand{\EE}{\end{equation}}

\theoremstyle{plain}
\newtheorem{theorem}{Theorem}

\newtheorem{lemma}[theorem]{Lemma}
\newtheorem{proposition}[theorem]{Proposition}

\theoremstyle{definition}
\newtheorem{definition}[theorem]{Definition}

\begin{document}
\title{Pattern-avoiding ascent sequences of length 3.}

\author{Andrew R. Conway\\Fairfield, Vic. 3078, Australia\\
\href{mailto:andrewascent@greatcactus.org}{\tt andrewascent@greatcactus.org},\\
Miles Conway\\Fairfield, Vic. 3078, Australia\\
\href{mailto: milesascent@greatcactus.org}{\tt milesascent@greatcactus.org}.\\
Andrew Elvey Price\\
CNRS, Institut Denis Poisson\\
Universit\'e de Tours, France\\
\href{mailto:andrew.elvey@univ-tours.fr}{\tt andrew.elvey@univ-tours.fr}.\\
Anthony J Guttmann\\
School of Mathematics and Statistics\\
The University of Melbourne\\
Vic. 3010, Australia\\
\href{mailto:guttmann@unimelb.edu.au}{\tt guttmann@unimelb.edu.au}
}

\date{}

\maketitle
\abstract
Pattern-avoiding ascent sequences have recently been related to set-partition problems and stack-sorting problems. While the generating functions for
several length-3 pattern-avoiding ascent sequences are known, those avoiding 000, 100, 110, 120 are not known. We have generated extensive series
expansions for these four cases, and analysed them in order to conjecture the asymptotic behaviour. 

We provide polynomial time algorithms for the $000$ and $110$ cases, and 
exponential time algorithms for the $100$ and $120$ case. We also describe
how the $000$ polynomial time algorithm was detected somewhat mechanically
given an exponential time algorithm.

For 120-avoiding ascent sequences we find that the generating
function has stretched-exponential behaviour and prove that the growth constant is the same as that for 201-avoiding ascent sequences, which is known. 

The other three generating functions have zero radius of convergence, which we also prove. For 000-avoiding ascent sequences we give what we believe to be the exact growth constant. We give the conjectured asymptotic behaviour for all four cases.

\section{Introduction}
\label{introduction}
 Given a sequence of non-negative integers, $n_1 n_2 n_3 \ldots n_k$  the number of {\em ascents} in this sequence is $$asc(n_1 n_2 n_3 \ldots n_k) = |\{ 1  \le j< k : n_j < n_{j+1} \} |.$$
  
  The given sequence is an {\em ascent sequence} of length $k$ if it satisfies $n_1=0$ and $n_i \in [0,1+asc(n_1 n_2 n_3 \ldots n_{i-1} )] $ for all $2 \le i \le k.$ For example, $(0,1,0,2,3,1,0,2)$ is an ascent sequence, 
  but $(0,1,2,2,4,3)$ is not, as $4 > asc(0122)+1=3.$ 
  
  Ascent sequences came to prominence when Bousquet-M\'elou et al. \cite{BCDK10} related them to $(2+2)$-free posets, and certain involutions, whose generating function was first given by Zagier \cite{Z01} . They have subsequently been linked to other combinatorial structures. See \cite{K11} for a number of examples. 
  The generating function for the number of ascent sequences of length $n$ is $$A(t) =\sum_{n \ge 0} a_n t^n =\sum_{n \ge 0} \prod_{i=1}^n (1-(1-t)^i)=1+t+2t^2+5t^3+15t^4+53t^5+217t^6+\cdots,$$
  and $$a_n \sim \frac{12\sqrt{3} \exp(\pi^2/12)}{\pi^{5/2}} n! \left ( \frac{6}{\pi^2} \right )^n \sqrt{n}.$$
  
  Later, Duncan and Steingrimsson \cite{DS11} studied {\em pattern-avoiding ascent sequences. }
  
  A pattern is simply a word on nonegative integers (repetitions allowed). Given an ascent sequence $(n_1 n_2 n_3 \ldots n_k)$,  a pattern $p$ is a subsequence $n_{i_1}n_{i_2}\ldots n_{i_j}$, where $j$ is just the length of $p$, and where the letters appear in the same relative order of size as those in $p.$ For example, the ascent sequence $(0,1,0,2,3,1)$ has three occurrences of the sequence $001,$ namely $002$, $003$ and $001$. If an ascent sequence does not contain a given pattern, it is said to be {\em pattern avoiding}. 
  
 The connection between pattern-avoiding ascent sequences and other combinatorial objects, such as set partitions, is the subject of \cite{DS11}, while the connection between pattern-avoiding ascent sequences and a number of stack sorting problems is explored in \cite{CCF20}.
 
 Considering patterns of length three, the number of ascent sequences of length $n$ avoiding the patterns $001,$ $010,$ $011,$ and $012$ is given in the OEIS \cite{OEIS} (sequence A000079) as $2^{n-1}.$ For the pattern $102$ the number is $(3^n+1)/2$ (OEIS A007051), while for $101$ and $021$ the number is just given by the $n^{th}$ Catalan number, $C_n,$ given in the OEIS as sequence A000108.
 
 More recently, the case of 210-avoiding ascent sequences, given in the OEIS as sequence A108304, was shown to be equivalent to the number of set partitions of $\{1,2,\cdots,n\}$ that avoid 3-crossings, the generating function for which was found by Bousquet-M\'elou and Xin \cite{BMX05}. It is D-finite, and the coefficients behave asymptotically as $$c210_n \sim \frac{3^9 \cdot 5}{2^5}\frac{\sqrt{3}}{\pi} \frac{9^n}{n^7}.$$
 
  More recently still, for the pattern 201 given in the OEIS as sequence A202062, Guttmann and Kotesovec \cite{GK21} found the generating function, which is not only D-finite but algebraic.
  The coefficients behave as $$c201_n \sim C\frac{\mu^n}{n^{9/2}},$$ where $\mu = 7.2958969432397\ldots$ is the largest root of the polynomial $x^3 - 8x^2 + 5x + 1,$ and the amplitude $C$ is
  $$C=\frac{35}{16}\left (\frac{4107}{\pi} - \frac{84}{\pi}\sqrt{9289}\cos\left (\frac{\pi}{3} +\frac{1}{3} \arccos\left [\frac{255709\sqrt{9289}}{24653006} \right ]\right )\right )^{1/2}.$$
 
This leaves the behaviour of just four remaining length-3 pattern-avoiding ascent sequences to be determined. They are 000, 100, 110 and 120. Quite short series for all four cases are given in the OEIS, but these are insufficient to conjecture the asymptotics.

We first developed an efficient dynamic programming algorithm to generate further terms, The efficiency of the algorithm is heavily pattern-dependent. For 110-avoiding ascent sequences we generated only 42 terms, but for 100-avoiding ascent sequences we generated 712 terms.

For 120-avoiding ascent sequences, we proved that the growth constant is the same as that for 201-avoiding ascent sequences, as given above. More generally, we found stretched-exponential behaviour, so that $$c120_n \sim C \cdot \mu^n \cdot \mu_1^{n^{3/8}} \cdot n^g,$$ where  $\mu = 7.2958969432397\ldots,$ $g=2,$ \ $\log{\mu_1}=-9.675 \pm 0.01$ and $C \approx 3700.$ The estimates of $C$ and $\mu_1$ depend sensitively on the validity of our estimate that $\sigma = 3/8$ exactly.

For 000-avoiding ascent sequences we found factorial behaviour, so that $$c000_n \sim C \cdot n!  \cdot \mu^n,$$ where $\mu \approx 0.27019$ and is conjectured to be $8/(3\pi^2)$ exactly.

For both 100-avoiding and 110-avoiding ascent sequences we found $$c_n \sim C \cdot \left ( \frac{3n}{4} \right )! \mu^n,$$ where both $C$ and $\mu$ are pattern dependent. We found that $\mu_{110} = 0.44 \pm 0.02,$  $\mu_{100} = 0.68 \pm 0.04,$ and more precisely that $\mu_{110}/\mu_{100} \approx 5/8.$

For these last three cases we give weak lower bounds that prove that factorial growth is to be expected in these cases. There are also presumably some sub-dominant terms, such as $n^g,$ but we were unable to estimate these.

In the next section we give details of our algorithm, and in the next four sections we study these four pattern-avoiding ascent sequences.

In the Appendix we describe the methods of series analysis used in this study.
 
\section{Sequence generation algorithm}

The sequences were generated by a set of slightly different dynamic programming algorithms. All of these are
restrictions added to a basic dynamic programming algorithm to enumerate the ascent sequences. The base algorithm will be explained
next, not because it is intrinsically useful but rather as the other algorithms are derived from it.

\subsection{Enumerating Ascent sequences}

\label{sec:basealg}

Consider a function $f(n,a,l)$ which 
gives the number of length $n$ suffixes of an ascent sequence where the
 prior portion of the sequence contains $a$ ascending pairs and the last number was $l$. This is useful as 
the number of ascent sequences of length $n$, $a_n=f(n-1,0,0)$, with $a_0=1$.

Consider all possibilities $i$ for the next number (the first number in the suffix), which must be between $0$
and $a+1$ inclusive. For each $i$, the rest of the suffix is of length $n-1$, and has a prior number $i$ and a prior number
of ascents of $a$ if $l\ge i$ and $a+1$ if $l<i$. This leads to a simple recursive definition:
$$
f(n,a,l) = \left\{ 
\begin{array}{ c l }
1 & \quad \textrm{if } n=0 \\
\sum_{i=0}^{a+1} f(n-1,a+\chi(l,i),i)   & \quad \textrm{otherwise}
\end{array}
\right.
$$
where
$$
\chi(l,i) = \left\{ 
\begin{array}{ c l }
1 & \quad \textrm{if } l<i \\
0 & \quad \textrm{otherwise.}
\end{array}
\right.
$$

This is trivial to implement in a recursive computer algorithm, and is very efficient using dynamic programming (storing each value, and not recomputing
any value already calculated). In particular, to enumerate $n$ terms, values of each of the
three arguments to $f$ never get above $n$, so the maximum number of terms visited is $O(n^3)$ and so the algorithm
uses time and space proportional to $O(n^3)$. This allows thousands of terms to be readily computed.

\begin{figure}
	\caption{Diagram of computation of $a_5=f(4,0,0)=53$ using the algorithm in section~\ref{sec:basealg}. 
		Each box represents a value of $f$ whose arguments are given in the top of the box. The value is in the lower left of the box.
	    Lines coming out of the right of a box go to the left of boxes corresponding to the instances of $f$ that are summed to produce the value in the box.
        Boxes and lines going to $n=0$ are omitted. }
    \label{fig:AscentDProg}
    \includegraphics[width=\textwidth]{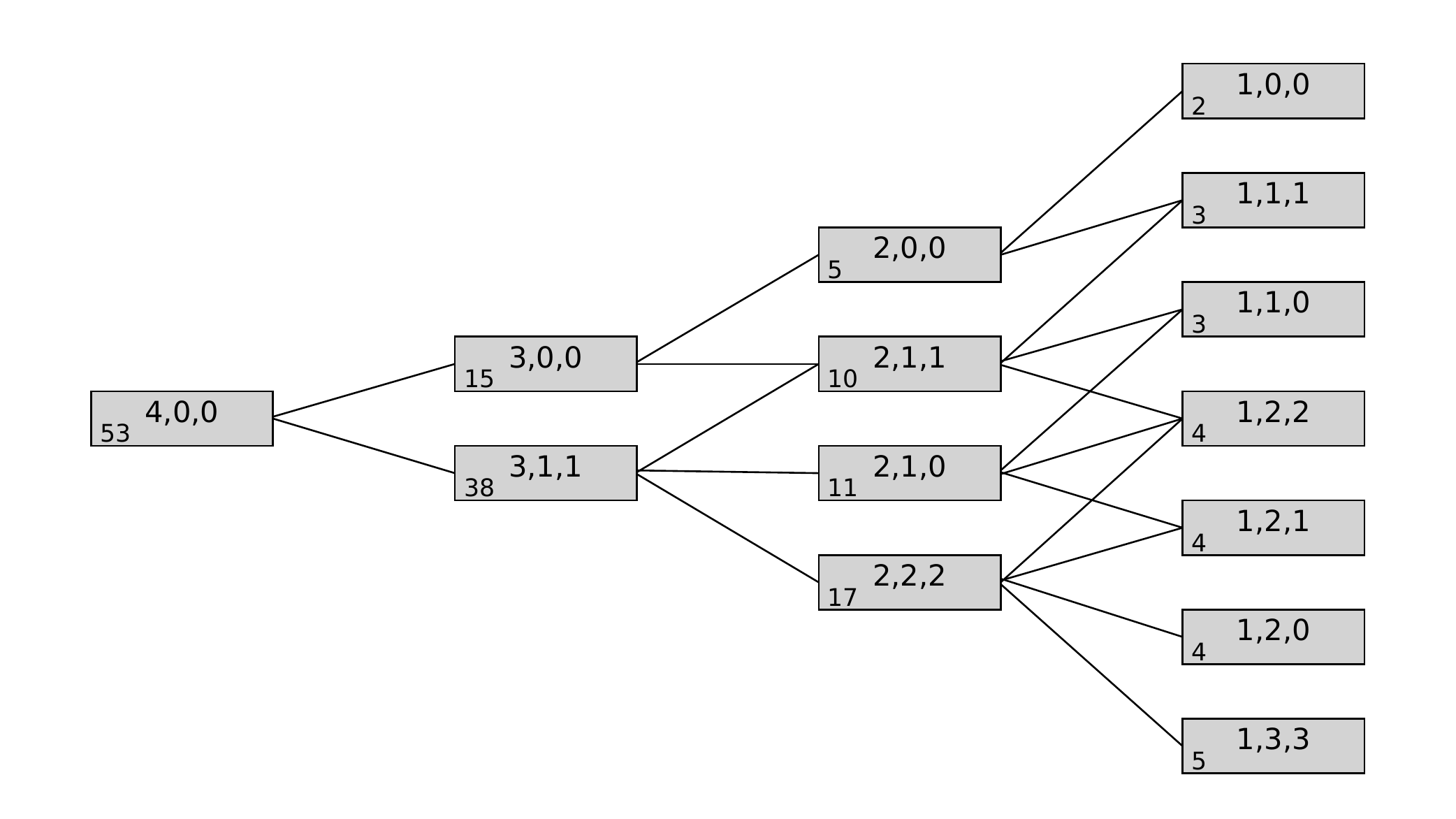}
\end{figure}

\begin{figure}
	\caption{Alternate diagram of computation of $a_5=f(4,0,0)=53$ using the algorithm in section~\ref{sec:basealg}. 
		Each box represents a value of $f$ whose arguments other than $n$ are given in the top of the box. The values for $n=0,1,\dots$ is in the lower left of the box.
		Lines coming out of the right of a box go to the left of boxes corresponding to the instances of $f$ that are summed to produce the value in the box. Each
        such line implicitly refers to one lower value of $n$. }
	\label{fig:AscentDProg2}
	\includegraphics[width=\textwidth]{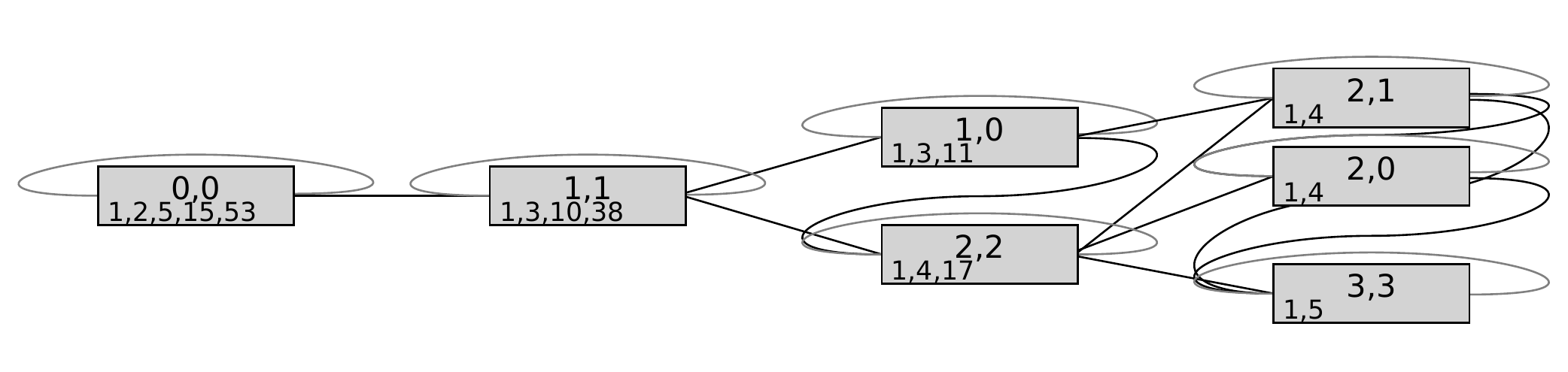}
\end{figure}

A diagram of the values actually computed to compute $a_5=53$ is given in figure~\ref{fig:AscentDProg}, and more
compactly in figure~\ref{fig:AscentDProg2} which uses the property of $f$ where the referenced elements are unchanged
(other than $n$ increasing) if $n$ is increased, except for the trivial case of $n=0$.

In future sections, to concentrate on the important part of the algorithm, the case $n=0$ will be omitted, and it
should always be assumed that the function will be $1$ in that case.

\subsection{Enumerating $000$ avoiding ascent sequences}

\label{sec:alg000}

A similar approach works for enumerating sequences avoiding the pattern $000$. Now when one considers adding a
number $i$ after some history, $i$ can be any value from $0$ to $a+1$ other than a value that
has been seen twice before. One could write a function that takes the same arguments as in
section~\ref{sec:basealg}, plus a set $S$ of the numbers seen exactly once before, and a set $P$ of proscribed numbers
present twice before. Then one could define a recursive function $g$
$$
g(n,a,l,S,P) = \sum_{i=0}^{a+1} \left\{ 
\begin{array}{ c l }
0 & \quad \textrm{if } i\in P \\
g(n-1,a+\chi(l,i),i,S\setminus\{i\},P\cup\{i\}) & \quad \textrm{if } i\in S \\
g(n-1,a+\chi(l,i),i,S\cup\{i\},P)   & \quad \textrm{otherwise.}
\end{array}
\right.
$$

This again could be computed in a straightforward manner. Both sets $S$ and $P$ can have $O(2^n)$ values, 
and so this is a much more computationally expensive algorithm. Fortunately there is a very straightforward
simplification. Any number in $P$ effectively does not exist, as far as the algorithm is concerned. 
However, the arguments are only used as numbers for their relative order and 0 element, not any other intrinsic 
numerical properties. It would be equally valid to rename the numbers such that the proscribed numbers
$P$ cease to exist and all other numbers are mapped to the integers starting from 0. That is, 
$$
g(n,a,l,S,P) = g(n,a-|P|,l-\sum_{i=0}^{l}\chi_P(i),S,\emptyset)
$$
where $\chi_P(i)$ is $1$ if $i\in P$ and $0$ otherwise.

This means there is no reason to remember the set $P$, allowing one to rewrite the recursive
equation in terms of a function $f^{000}(n,a,l,S)=g(n,a,l,S,\emptyset)$:
$$
f^{000}(n,a,l,S) = \sum_{i=0}^{a+1} \left\{ 
\begin{array}{ c l }
f^{000}\left(n-1,a+\chi(l,i)-1,i-1,r(S\setminus\{i\},i)\right) & \quad \textrm{if } i\in S \\
f^{000}(n-1,a+\chi(l,i),i,S\cup\{i\})   & \quad \textrm{otherwise}
\end{array}
\right.
$$
where $r(S,i)$ is the renumbering function that takes a set $S$, and reduces 
the value of each element greater than $i$ by 1, as the old number $i$ is edited out of
existence.

Note that this erasure of numbers out of existence may mean that $a$ or $l$ ends up being $-1$, which doesn't
need any special handling; $a=-1$ means the next value must be a $0$, and $l=-1$ means that the next value {\it will} be larger than $l$.

Logistically, the set $S$ is represented on a computer as a long integer, where the $i^th$ bit
is $1$ iff $i\in S$. Then $r(S,i)$ can be easily done via bit masking and shifting.

The desired sequence $\{c000_n \}=\{1,2,4,10,27,83,277,\ldots\}$ is then given by $c^{000}_n=f^{000}(n-1,0,0,\{0\})$.

\begin{figure}
	\caption{Diagram of the first 6 layers (sufficient to compute 6 terms) of the computation of $000$ avoiding ascent sequences using the algorithm in section~\ref{sec:alg000}. 
		The set $S$ is represented as a binary string in the third argument at the top of each box; $i\in S$ iff the $i^{th}$ digit from the right is a 1. 
		Otherwise the interpretation is the same as figure~\ref{fig:AscentDProg2}. }
	\label{fig:Avoid000DP}
	\includegraphics[width=\textwidth]{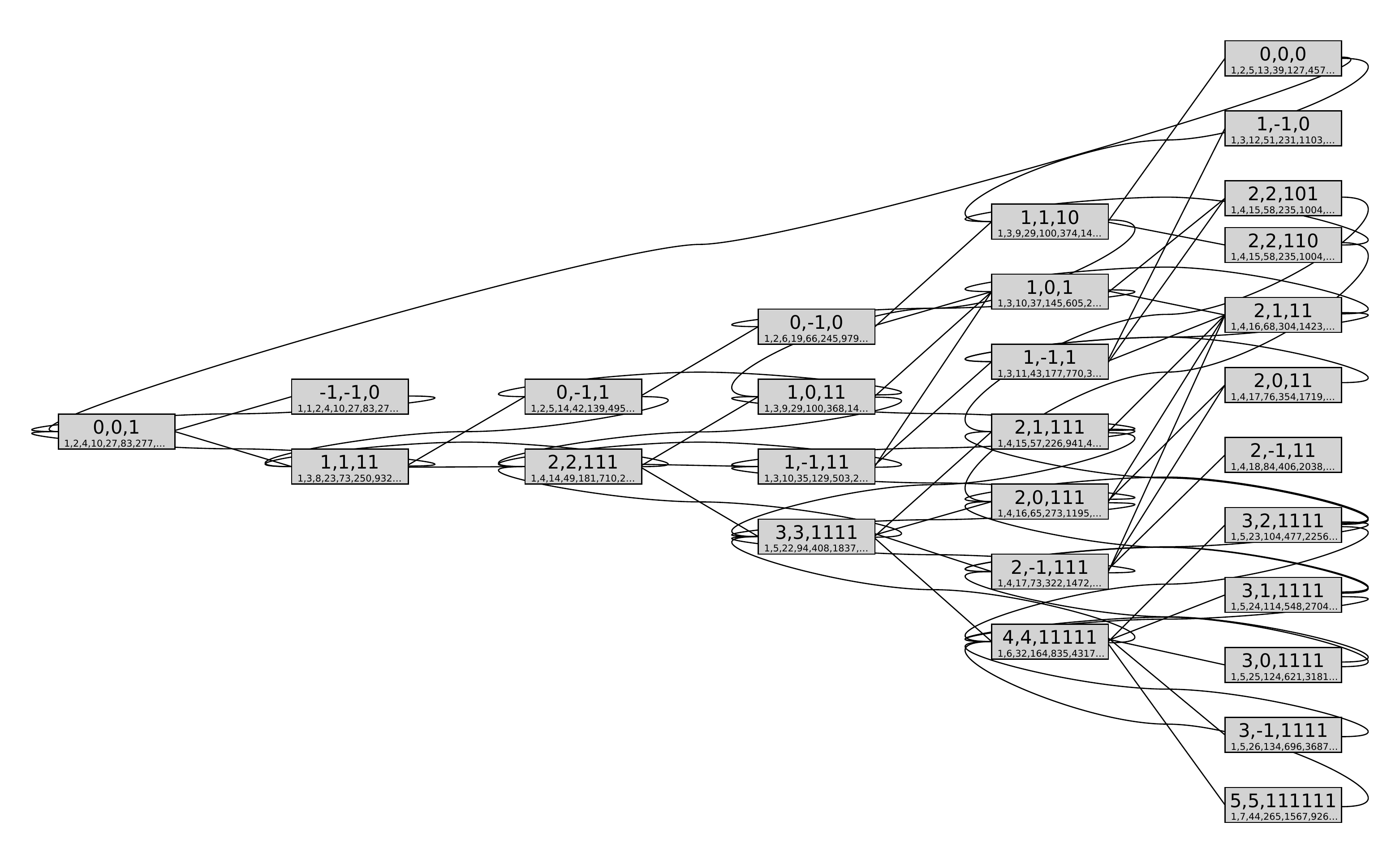}
\end{figure}

The number of possible values of $S$ is no more than $2^n$, so the algorithm is no worse
than $O(n^3 2^n)$ and in practice is slightly better. It can be readily calculated by this algorithm to about 30 terms. 

A graph of the structure of this computation is given in figure~\ref{fig:Avoid000DP}. Note that the
values corresponding to $a=2,l=2,S={0,2}$ are the same as $a=2,l=2,S={1,2}$ for all $n$. Similar
behaviour happens more frequently as more levels are shown; indeed it turns out empirically
that $f^{000}(n,a,l,S)$ for given values of $n$, $a$, and $l$ only depends upon the cardinality of $S$.
This implies a much more efficient algorithm yet - see section~\ref{sec:better000}.

\subsection{Enumerating $100$ avoiding ascent sequences}

\label{sec:alg100}

A somewhat similar approach works for avoiding $100$. In this case we want to rewrite a number
out of existence if a number is ever encountered that is lower than any previously seen number.
This can be done by keeping track of the largest number so far seen, $m$. Then

$$
f^{100}(n,a,l,m) = \sum_{i=0}^{a+1} \left\{ 
\begin{array}{ c l }
f^{100}\left(n-1,a+\chi(l,i)-1,i-1,m-1\right) & \quad \textrm{if } i<m \\
f^{100}(n-1,a+\chi(l,i),i,i)   & \quad \textrm{otherwise.}
\end{array}
\right.
$$

The desired sequence is then given by $c^{100}_n=f^{100}(n-1,0,0,0)$.

This directly provides an efficient algorithm, $O(n^4)$, which enables many hundreds of terms
to be readily computed.

\subsection{Enumerating $110$ avoiding ascent sequences}

\label{sec:alg110}

The $110$ avoiding case is very similar to the $000$ avoiding case, except this time
when we get a repeated number, we want to write out of existence any number less than
it.

$$
f^{110}(n,a,l,S) = \sum_{i=0}^{a+1} \left\{ 
\begin{array}{ c l }
f^{110}\left(n-1,a+\chi(l,i)-i,0,r(S,i)\right) & \quad \textrm{if } i\in S \\
f^{110}(n-1,a+\chi(l,i),i,S\cup\{i\})   & \quad \textrm{otherwise.}
\end{array}
\right.
$$
In this case the renumbering function $r(S,i)$ removes any value in $S$ less than $i$, and reduces the values of all others by $i$.

The desired sequence is then given by $c^{110}_n=f^{110}(n-1,0,0,\{0\})$.

Implementing this as a dynamic programming algorithm has the same upper bound as the $000$ avoiding algorithm described in section \ref{sec:alg000},
although in practice we write out of existence more numbers and fewer states occur in practice, allowing 40 or so terms
to be readily computed with current hardware.

\subsection{Enumerating $120$ avoiding ascent sequences}

\label{sec:alg120}

The $120$ avoiding case is very similar to the $110$ avoiding case,
except when we encounter a value $i$ larger than some previously seen value,
we want to erase out of existence all numbers smaller than the largest previously
seen value less than $i$.

$$
f^{120}(n,a,l,S) = \sum_{i=0}^{a+1} 
f^{120}\left(n-1,a+\chi(l,i)-s(S,i),i-s(S,i),r(S,s(S,i))\right).
$$
In this case the renumbering function $r(S,i)$ has the same meaning as in section~\ref{sec:alg110} and the 
function $s(S,i)$ means the largest element in $S$ smaller than $i$, or $0$ if there is none. 
Note that in practice $S$ will always contain the element $0$. 

The desired sequence is then given by $c^{120}_n=f^{120}(n-1,0,0,\{0\})$.

It is more difficult to get large values of $a$ in this case than prior cases, as two consecutive
increases of number in the sequence (increasing $a$ by $2$) will cause the first to be rewritten out of existence, reducing
$a$ by 1. This means $n$ must be increased by 2 to increase the maximum value of $a$ by 1. 
This is primarily important as the maximum element in $S$ is determined by the maximum value of $a$,
so the algorithm becomes $O(n^3 2^{n/2})$. This allows about twice as many terms as the 000 or 110 
algorithms, or in the seventies in practice with current hardware.

\subsection{Better $000$ algorithm}
\label{sec:better000}

This section presents a more efficient algorithm for the $000$ case
than presented in section~\ref{sec:alg000}. 
Perhaps more interesting than the algorithm itself is the method
used to discover it.

For many years, we have suspected that looking for frequently repeated large numbers
in the dynamic programming cache will lead to the observation that a similar, more efficient version
of the same algorithm exists, tracking a subset of the state information
that was thought to be needed. This is the first time we have actually seen strong evidence of this, with the majority of large numbers repeated.

Extensive numerical evidence demonstrated that, in the $000$ algorithm presented above, the value of 
$f^{000}(n,a,l,S)$ for given values of $n$, $a$, and $l$ is the same
for many different values of $S$ with the same cardinality.

To see why, consider a set $T$ and non-negative integer $i$ such that $i\notin T$ and $i+1\notin T$.
We will demonstrate that $f^{000}(n,a,l,T\cup\{i\})=f^{000}(n,a,l,T\cup\{i+1\})$
if $i<l$. Consider a specific suffix $u$ counted by $f^{000}(n,a,l,T\cup\{i\})$. Find each maximal contiguous subsequence in $u$
containing just $i$ and $i+1$. Reverse each of these sequences, and replace
each $i$ by $i+1$. The resulting suffix $v$ is counted by $f^{000}(n,a,l,T\cup\{i+1\})$
as the number of $i$ and $i+1$ values are swapped, other values are unchanged; the number of 
ascents is unchanged (including at the start as neither $i$ or $i+1$
can exceed $l$), and $i+1$ will be allowed as $l$ was allowed and $i+1\le l$.
Furthermore, this is a bijection, so $f^{000}(n,a,l,T\cup\{i\})=f^{000}(n,a,l,T\cup\{i+1\})$.

This can be used to canonicalise the value of $S$ used in the 
recursive definition of $f^{000}$, decreasing the number
of states visited. In particular, for the
case of a duplication where a number is rewritten out of existence,
remove that number as usual. When a new value $i$ is added, then 
instead of evaluating for $S+{i}$, bubble the value $i$ down using
multiple invocations of the prior paragraph until the value below it
is already in $S$. Define $S$ to be \emph{compacted} if it is a (possibly empty) set
of consecutive integers starting at zero. In both cases, assuming the S input to $f^{000}$ is compacted, then all calls to $f^{000}$ it produces
will also be compacted. As the initial call to $f$ is the compacted set $\{0\}$, all calls will be of compacted sets.

A compacted set can be represented by its cardinality, which will not
exceed $n$ for enumerating $n$ terms. This means the enumeration
algorithm becomes $O(n^4)$ which is much more efficient, and allows
easy enumeration of hundreds of terms on current hardware.

This is useful as it enables us to generate vastly more terms of
the series; it is also interesting as it demonstrates how a mechanical
operation (checking the dynamic programming cache for repeated large
values) can lead to a polynomial time algorithm given an exponential
algorithm. A mechanical method of getting good ideas... or at least
becoming aware of their existence, is of great value.

The $110$ algorithm has very few repeated large numbers. The $120$
algorithm has an intermediate amount of repeated large numbers. For instance,
$f^{120}(-,4,0,{0,1,2,4})=f^{120}(-,4,0,{0,1,3,4})=f^{120}(-,4,0,{0,2,3,4})$
for all values of $n$ tried (going $1,6,32,160,778,3747\cdots$ with $n$). However the relationship is more complex than the $000$ case, the efficiency gains are lower,
and we can already enumerate many terms for this sequence anyway,
so we did not pursue a new $120$ algorithm.

In the next four sections, we analyse the extended sequences produced by these algorithms, in order to conjecture the asymptotic behaviour, in each case.

\section{120-avoiding ascent sequence}\label{120}
This sequence is given as A202061 in the OEIS \cite{OEIS} to order O$(x^{14})$, and we have extended this to O$(x^{73}).$ We have used these exact coefficients to derive 200 further coefficients by the method of series extension \cite{G16}, and briefly described in the Appendix. 

We first plot the ratios of the coefficients $r_n = c_n/c_{n-1}$ against $1/n.$ If one has a pure power-law type singularity, such a plot should be linear, with ordinate interception giving an estimate of the growth constant. 

\begin{figure}[h!] 
\begin{minipage}[t]{0.45\textwidth} 
\centerline{\includegraphics[width=\textwidth]{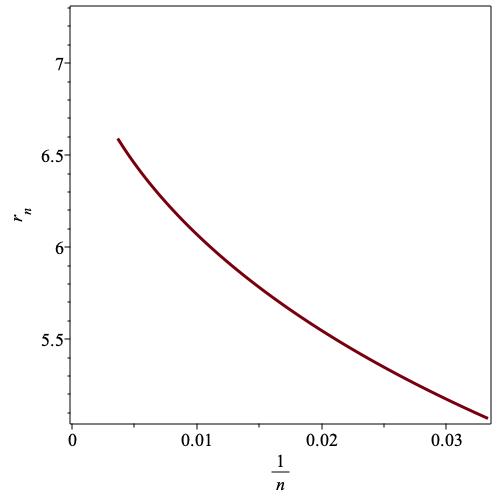}}
%\captionsetup{justification=centering}
\caption{
Plot of ratios against $1/n$ for $120$-avoiding ascent sequences.} 
\label{fig:rat1}
\end{minipage}
\hspace{0.05\textwidth}
\begin{minipage}[t]{0.45\textwidth} 
\centerline {\includegraphics[width=\textwidth]{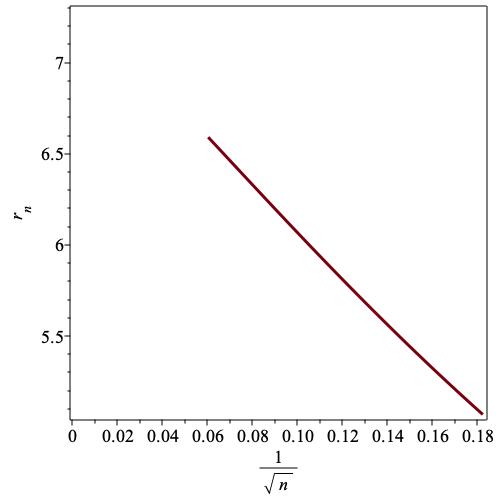}}
% \captionsetup{justification=centering}
\caption{Plot of ratios against $1/\sqrt{n}$ for  $120$-avoiding ascent sequences. }
\label{fig:rat2}
\end{minipage}
\end{figure}

The ratio plot is shown in Fig. \ref{fig:rat1}, and it displays considerable curvature. By contrast, if the ratios are plotted against $1/\sqrt{n},$ as shown in  Fig. \ref{fig:rat2}, the plot is virtually linear, and intercepts the ordinate at around 7.3, which is our initial estimate of the growth constant.

This behaviour of the ratios suggests that the singularity is of stretched-exponential type, so that the coefficients behave as 
\BE \label{eq:ana}
c_n \sim C \cdot \mu^n \cdot \mu_1^{n^\sigma} \cdot n^g,
\EE
 with $\mu \approx 7.3 $ and $\sigma \approx 1/2.$ Given such a singularity, the ratios will behave as
\begin{multline} \label{eq:rn}
r_n = \mu \left (1 + \frac{\sigma \log \mu_1}{n^{1-\sigma}} + \frac{g}{n} + \frac{\sigma^2 \log^2 \mu_1}{2n^{2-2\sigma}} + \frac {(\sigma-\sigma^2)\log \mu_1+2g\sigma \log \mu_1}{2n^{2-\sigma}} \right . \\
 \left . {}+ \frac{\sigma^3 \log^3 \mu_1}{6n^{3-3\sigma}} +{\rm O}(n^{2\sigma-3}) + {\rm O}(n^{-2}) \right ).
%\notag
\end{multline}

\begin{figure}[h!] 
\begin{minipage}[t]{0.45\textwidth} 
\centerline{
\includegraphics[width=\textwidth]{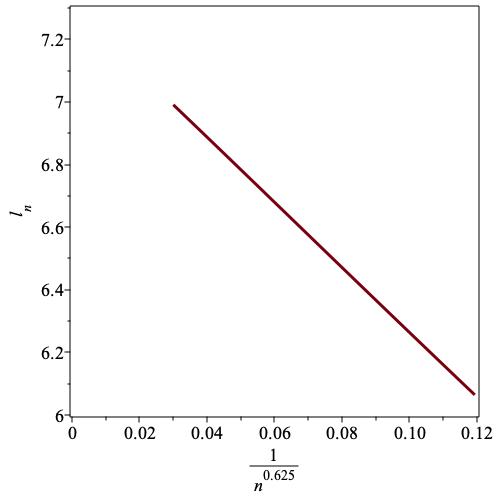}}
%\captionsetup{justification=centering}
\caption{Plot of linear intercepts $l_n$ against $1/n^{5/8}$ for $120$-avoiding ascent sequences.}
 \label{fig:ln}
\end{minipage}
\hspace{0.05\textwidth}
\begin{minipage}[t]{0.45\textwidth} 
\centerline {\includegraphics[width=\textwidth]{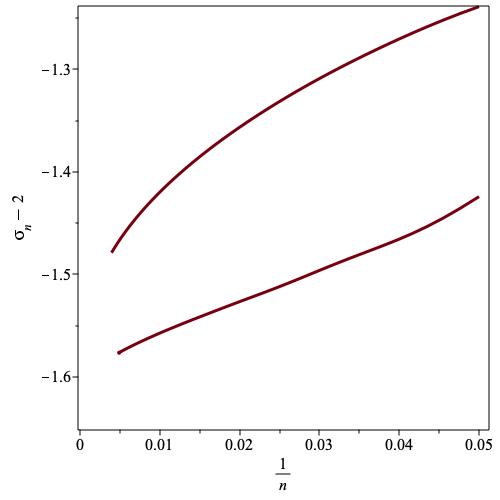}}
% \captionsetup{justification=centering}
\caption{Plot of estimators of $2-\sigma$ against $1/n$ for $120$-avoiding ascent sequences.} 
\label{fig:ln2}
\end{minipage}
\end{figure}

One can eliminate the O$(1/n)$ term in the expression for the ratios by studying instead the linear intercepts, $$l_n \equiv n\cdot r_n - (n-1)\cdot r_{n-1} \sim \mu \left (1 + \frac{\sigma^2 \log \mu_1}{n^{1-\sigma}} + \frac{\sigma^2(2\sigma-1) \log^2 \mu_1}{2n^{2-2\sigma}} \right ).$$ When we plot $l_n$ against $1/\sqrt{n}$ there is still some curvature in the plot, butthis disappears when we plot $l_n$ against $1/n^{5/8},$ as shown in Fig. \ref{fig:ln}. This suggests that $\sigma \approx 3/8$ is closer than $1/2.$

In order to estimate $\sigma$ without knowing or assuming $\mu,$ we can use one (or both) of the following estimators:

From eqn. (\ref{eq:rn}), it follows that 
\BE \label{eq:sig1a}
r_{\sigma_n} \equiv \frac{r_n}{r_{n-1}} \sim 1 + \frac{(\sigma-1)\log{\mu_1}}{n^{2-\sigma}} + {\rm O}(1/n^2),
\EE
so $\sigma$ can be estimated from a plot of $\log(r_{\sigma_n}-1)$ against $\log{n},$ which should have gradient $\sigma-2.$ The local gradients can be calculated and plotted against $1/n^\sigma,$ using any approximate value of $\sigma.$

Another estimator of $\sigma$ when $\mu$ is not known follows from eqn. (\ref{eq:ana}):
\BE \label{eq:sig2a}
a_{\sigma_n} \equiv \frac{c_n^{1/n}}{c_{n-1}^{1/(n-1)}} \sim 1 + \frac{(\sigma-1)\log{\mu_1}}{n^{2-\sigma}} + {\rm O}(1/n^2),
\EE
so again $\sigma$ can be estimated from a plot of $\log(a_{\sigma_n}-1)$ against $\log{n}.$ Again, estimates of $\sigma$ are found by extrapolating the local gradient  against $1/n^\sigma.$

%\centerline{\includegraphics[width=\textwidth,angle=0]{Av25314rat.jpg} }

While these two estimators are equal to leading order, they differ in their higher-order terms. We show these two estimators in Fig. \ref{fig:ln2}), and both estimators are consistent with the estimate $2-\sigma \approx -1.625,$ so that $\sigma \approx 3/8.$
We cannot of course exclude nearby values, such as $\sigma = 0.4.$

If $\sigma$ is known, or assumed, one can estimate at least some of the critical parameters by direct fitting. In particular, one can fit the ratios to
\BE \label{eqn:ratfit}
r_n = c_1+\frac{c_2}{n^{1-\sigma}} +\frac{c_3}{n}+\frac{c_4}{n^{2-2\sigma}},
\EE
by solving the linear system obtained by taking four successive ratios $r_{k-2},\,r_{k-1},\,r_k,\,r_{k+1},$ from which one can estimate the parameters $c_1, \cdots,c_4.$ One increases $k$ until one runs out of known (or estimated) coefficients.
Then $c_1$ estimates $\mu,$ $c_2$ estimates $\mu\cdot \sigma \log(\mu_1)$, $c_3$ estimates $\mu\cdot g$ (assuming $\sigma\ne 1/2$),  and $c_4$ gives estimators of $\mu\cdot \sigma^2 \log^2(\mu_1)/2.$ 

\begin{figure}[h!] 
\begin{minipage}[t]{0.45\textwidth} 
\centerline{\includegraphics[width=\textwidth]{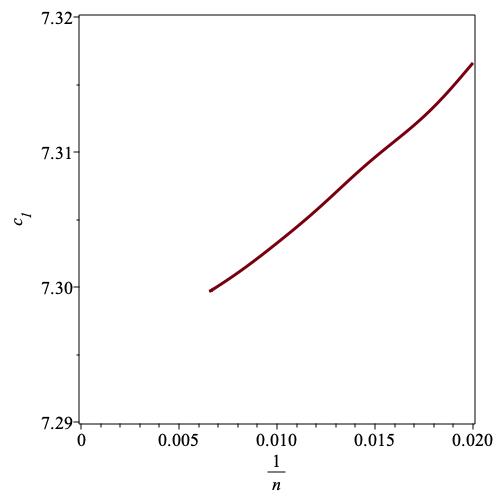}}
%\captionsetup{justification=centering}
\caption{
Plot of estimates of $c_1$ against $1/n$ for $120$-avoiding ascent sequences.} 
\label{fig:c1}
\end{minipage}
\hspace{0.05\textwidth}
\begin{minipage}[t]{0.45\textwidth} 
\centerline {\includegraphics[width=\textwidth]{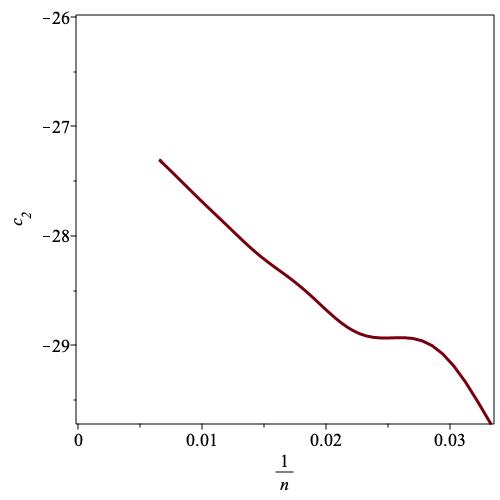}}
% \captionsetup{justification=centering}
\caption{Plot of estimates of $c_2$ against $1/n$ for $120$-avoiding ascent sequences.}
\label{fig:c2}
\end{minipage}
\end{figure}

\begin{figure}[h!] 
\begin{minipage}[t]{0.45\textwidth} 
\centerline{
\includegraphics[width=\textwidth]{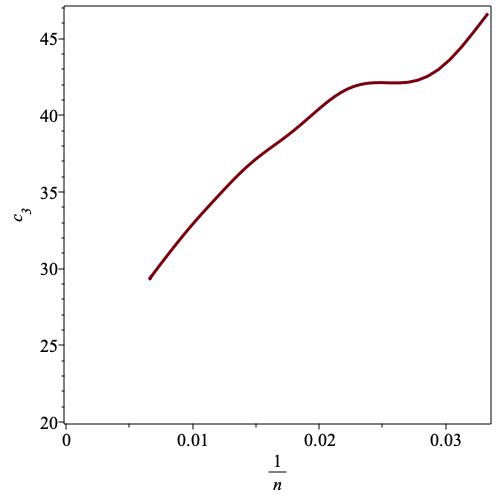}}
%\captionsetup{justification=centering}
\caption{Plot of estimates of $c_3$ against $1/n$ for $120$-avoiding ascent sequences..}
 \label{fig:c3}
\end{minipage}
\hspace{0.05\textwidth}
\begin{minipage}[t]{0.45\textwidth} 
\centerline {\includegraphics[width=\textwidth]{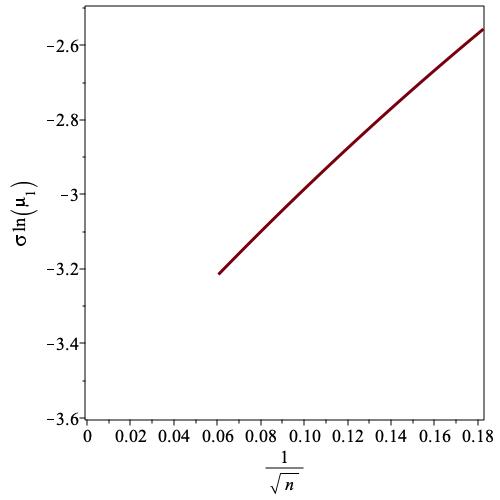}}
% \captionsetup{justification=centering}
\caption{Plot of estimators of $\sigma \cdot \log{\mu_1}$ against $1/\sqrt{n}$ for $120$-avoiding ascent sequences.} 
\label{fig:mu1}
\end{minipage}
\end{figure}

We show the results of doing this, assuming $\sigma=3/8,$ in Figs \ref{fig:c1}, \ref{fig:c2}, \ref{fig:c3}, from which we estimate $c_1 \approx 7.295,$ $c_2 \approx -26.5,$ and $c_3 \approx -20.$ 

From these it follows that $\mu \approx 7.295, $
$\log{\mu_1} \approx -9.68,$ and $g \approx 2.7.$ If we repeat this analysis assuming $\sigma = 0.4,$ the estimate of $c_1$ barely changes, increasing to 7.297, but $c_2 \approx -20.5,$ and $c_3 < 0,$ so $\log{\mu_1} \approx  -7,$ and $g < 0.$ So these last two parameters are seen to be very sensitive to the assumed value of $\sigma.$

From eqn. (\ref{eq:rn}), if we know (or conjecture) $\mu$ and $\sigma,$ we can use this to estimate $\mu_1,$ as
\BE \label{eq:mu1}
\left ( \frac{r_n}{\mu} - 1\right )\cdot n^{1-\sigma} \sim \sigma \cdot \log(\mu_1).
\EE
In Fig. \ref{fig:mu1} we show the relevant plot, where we have used the estimates $\mu=7.295$ and $\sigma=3/8.$ We estimate the ordinate intercept to be at around -3.6, from which follows $\log{\mu_1} \approx -9.6.$
This is in agreement with the estimate obtained above by direct fitting with $\sigma=3/8$ assumed. Similarly, if we assume $\sigma=0.4,$ we find $\log{\mu_1} \approx -7.2,$ again in agreement with the value found above by direct fitting.

%\centerline{\includegraphics[width=\textwidth,angle=0]{Av25314rat.jpg} }

Recall that the growth constant for ascent sequences is $6/\pi^2.$ It is interesting to note that $72/\pi^2=7.2951\ldots.$ However the growth constant of 201-avoiding ascent sequences \cite{GK21} is $7.2958969\ldots ,$ so either value of the growth constant appears possible. However, we can prove that 201-ascent sequences have the same growth constant as 120-avoiding ascent sequences, as shown in Sec. \ref{AEP}. 

Using this knowledge, we can repeat the above analysis to establish the value of the exponent $\sigma,$ but incorporating the known value of $\mu.$ Doing this, we find $\sigma = 0.372 \pm 0.006,$ which accords with our conjecture above that $\sigma=3/8$ exactly.

%We can also now compare this sequence with the behaviour of the coefficients of ascent sequences, after dividing the ascent-sequence coefficients by $n!$, and constructing the Hadamard quotient of the coefficients of the two sequences. That is to say, the ascent-sequence coefficients are known to grow as $a_n/n! \sim C  \cdot (6/\pi^2)^n \cdot \sqrt{n},$ and our best estimate of the behaviour of the coefficients of $120$-avoiding ascent sequences is $c_n \sim D \cdot n! \cdot \mu_1^{n^\sigma}\cdot  \mu^n \cdot n^g.$ The Hadamard quotient  is$$h_n \equiv \frac{c_n n!}{a_n} \sim const. \lambda^n \cdot \mu_1^{n^\sigma}  \cdot n^{g-1/2},$$ where $\lambda = \mu \pi^2/6.$ Extrapolation of the ratios $h_n/h_{n-1}$ should give estimates of $\lambda.$ We have eliminated terms of O$(1/n)$ and O$(1/n^2)$ and then plot linear intercepts against $1/n^{1-\sigma}$ in Fig. \ref{fig:hrat120}. Doing this gives $\lambda \approx 12.000,$ which implies $\mu = 72/\pi^2,$ as suggested earlier, though with greater precision than previous estimates.

\begin{figure}[h!] 
\begin{minipage}[t]{0.45\textwidth} 
\centerline{
\includegraphics[width=\textwidth]{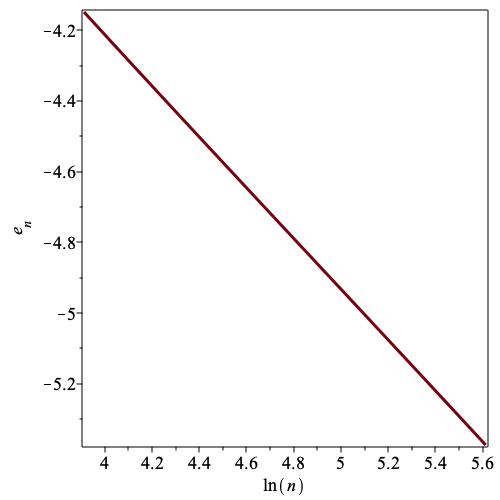}}
%\captionsetup{justification=centering}
\caption{Plot of $e_n$ against $\log{n}$ for $120$-avoiding ascent sequences.}
 \label{fig:e3}
\end{minipage}
\hspace{0.05\textwidth}
\begin{minipage}[t]{0.45\textwidth} 
\centerline {\includegraphics[width=\textwidth]{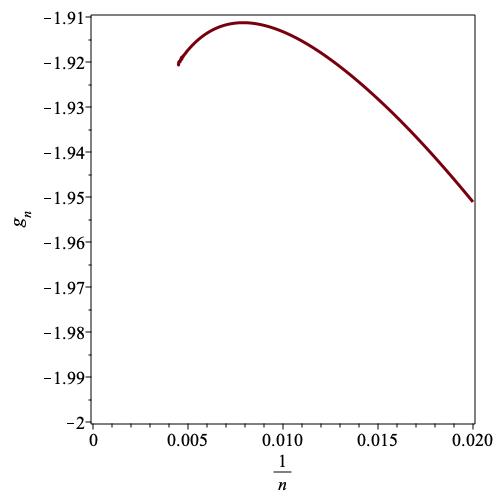}}
% \captionsetup{justification=centering}
\caption{Plot of estimators of exponent $-g$ against $1/n$ for $120$-avoiding ascent sequences.} 
\label{fig:g1}
\end{minipage}
\end{figure}

Having established the value of $\mu$ and conjectured the value of $\sigma,$ we are now in a better position to estimate the other parameters. We define the normalised coefficients $d_n \equiv c_n/\mu^n.$ Then we have
$$\log{d_n} \sim \log{C} + n^\sigma \cdot \log{\mu_1} +g\log{n},$$
$$e_n \equiv ((n-1)^\sigma \log{d_n} -n^\sigma \log{d_{n-1}})\cdot n^{1-\sigma} \sim -\sigma \log{C} + g - \sigma \cdot g \log{n}.$$
So a plot of $e_n/\sigma$ vs. $\log{n}$ should be linear, with gradient $-g.$ We show this plot in Fig. \ref{fig:e3}, and plot the local gradients against $1/n$ in Fig. \ref{fig:g1}. We conclude from this that $g \approx 2.$

Using this value of $g,$ we can get a more precise estimate of $\mu_1.$ Define the newly normalised coefficients $$f_n \equiv c_n/(n^g \mu^n).$$ Then
$$\frac{\log{f_n} -\log{f_{n-1}}}{n^\sigma-(n-1)^\sigma} \sim \log{\mu_1}.$$

\begin{figure}[h!] 
\
\centerline{
\includegraphics[width=0.5\textwidth]{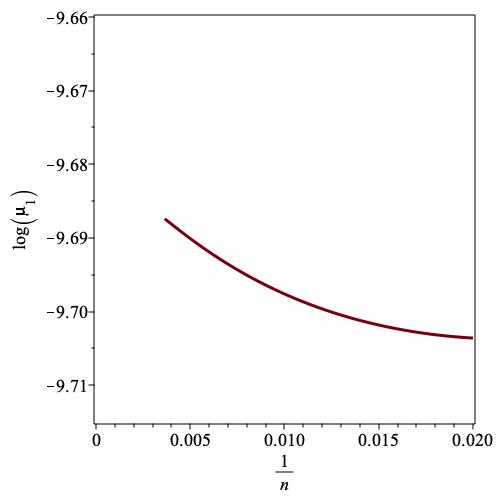}}
%\captionsetup{justification=centering}
\caption{Plot of $\log{\mu_1}$ against $1/n$ for $120$-avoiding ascent sequences.}
 \label{fig:lmu1}
\end{figure}

We show in Fig. \ref{fig:lmu1} a plot of these estimates of $\log{\mu_1}$ against $1/n,$ and from this we estimate $\log{\mu_1} = -9.675 \pm 0.01.$

Finally, we estimate the constant $C$ by dividing the coefficients by $\mu^n \cdot \mu_1^{n^\sigma} \cdot n^g$ with the assumed values of the parameters $\mu,\,\, \sigma,\,\,$ and $\mu_1.$ In this way we estimate $C \approx 3700.$

We conclude this section by giving our best estimate for the asymptotic behaviour of  the coefficients of $120$-avoiding ascent sequences as $$c_n \sim C \cdot \mu^n \cdot \mu_1^{n^\sigma} \cdot n^g,$$ where $\mu = 7.2958969\ldots,$ and is the largest solution of the cubic equation $x^3-8x^2+5x+1,$ $\sigma \approx 3/8,$ $\log{\mu_1} = -9.675 \pm 0.01,$  $g =2 $ and $C \approx 3700.$ Apart from the value of the growth constant $\mu,$ the other parameters, $\mu_1$ and $C$ depend sensitively on the precision of our estimate of $\sigma.$ If $\sigma \ne 3/8$ our estimates of $\mu_1$ and $C$ should not be believed.

\section{Proof that 120-avoiding and 201-avoiding ascent sequences have the same growth constant} \label{AEP}
In this section we show that for certain patterns $\pi$ the counting sequence of $\pi$-avoiding ascent sequences has the same growth rate as the counting sequence of simpler objects, which we call {$\pi$-avoiding weak ascent sequences}. Moreover, weak ascent sequences have a symmetry property which allows us to show that 120-avoiding and 201-avoiding ascent sequences have the same growth constant.

\begin{definition}
A sequence of non-negative integers, $n_1 n_2 n_3 \ldots n_k$ is a {\em weak ascent sequence} of length $k$ if it satisfies 
\[\max\{n_1,n_2,\ldots,n_k\} \leq asc(n_1 n_2 n_3 \ldots n_{k} ).\]
\end{definition}

\begin{lemma}If $n_{1}n_{2}\ldots n_{k}$ is an ascent sequence, then it is also a weak ascent sequence.
\end{lemma}
\begin{proof}
Assume $n_{1}n_{2}\ldots n_{k}$ is an ascent sequence and let $i\in[1,n]$ satisfy $n_{i}=\max(\{n_1,n_2,\ldots,n_k\})$ with $i$ minimal. If $i=1$, then $n_{i}=n_{1}=0$, so $n_{j}=0$ for all $j\in[1,n]$ and so the sequence is a weak ascent sequence.

If $i>1$, then by the definition of {\em ascent sequence}, we have 
\[n_{i}\le  asc(n_1 n_2 n_3 \ldots n_{i-1} )+1.\]
Moreover, by the definition of $i$, we have $n_{i-1}<n_{i}$, so $asc(n_1 n_2 n_3 \ldots n_{i-1} )+1=asc(n_1 n_2 n_3 \ldots n_{i} )$. Combining these yields
\[\max\{n_1,n_2,\ldots,n_k\}=n_{i}\le asc(n_1 n_2 \ldots n_{i} )\leq asc(n_1 n_2 \ldots n_{k} ),\]
so $n_{1}n_{2}\ldots n_{k}$ is a weak ascent sequence.
\end{proof}
\begin{proposition} \label{prop:weak-ascent}
Let $j\geq3$ and let $\pi=\pi_1\pi_2\ldots\pi_{j}$ be a sum-indecomposable permutation of $01\ldots j-1$,  that is, there is no $i\in[1,j-1]$ satisfying $\max\{\pi_1,\ldots,\pi_i\}<\min\{\pi_{i+1},\ldots,\pi_j\}$. Let $w_{k}$ be the number of $\pi$-avoiding weak ascent sequences and let $c_{k}$ be the number of $\pi$-avoiding ascent sequences. Then the exponential growth rates $\mu_{c}=\lim_{k\to\infty}\sqrt[k]{c_{k}}$ and $\mu_{w}=\lim_{k\to\infty}\sqrt[k]{w_{k}}$ exist and are equal to each other.
\end{proposition}
\begin{proof}
We start by showing that the limits exist, by proving that $c_{m+n}\geq c_{m}c_{n}$ and then applying Fekete's Lemma. If $C_{1}$ and $C_{2}$ are $\pi$-avoiding ascents sequences with lengths $m$ and $n$ and maximum values $N_{1}$ and $N_{2}$ respectively, then we define $\tilde{C}_{2}=C_{2}+N_{1}$ to be the sequence defined by increasing each value of $C_{2}$ by $N_{1}$. Then $C_{1}$ must contain at least $N_{1}$ ascents, so $C=C_{1}\tilde{C}_{2}$ is an ascent sequence, and it is $\pi$-avoiding because $\pi$ is sum-indecomposable. Since each pair $C_{1}$, $C_{2}$ defines a distinct sequence $C$ of length $m+n$, this implies that $c_{m+n}\geq c_{m}c_{n}$. By exactly the same argument, $w_{m+n}\geq w_{m}w_{n}$. Now Fekete's lemma implies that the limits 
\[\mu_{c}=\lim_{k\to\infty}\sqrt[k]{c_{k}}~~~~~\text{and}~~~~~\mu_{w}=\lim_{k\to\infty}\sqrt[k]{w_{k}}\] exist, although we do not prove that they are necessarily finite.
 
We will now show that
\[c_{k^2+2k}\geq w_{k}^{k},\]
as then combining this with $c_{k}\leq w_{k}$ and taking limits yields the desired result.

Let $W_{1},W_{2},\ldots,W_{k}$ be $\pi$-avoiding weak ascent sequences of length $k$, with maximum values $N_{1},N_{2},\ldots,N_{k}$, respectively. There are $w_{n}^{n}$ choices of these, so we just need to show that we can construct a unique $\pi$-avoiding ascent sequence $C$ of length $n^2+n$ with each such choice. We construct $C$ as follows:
\[C=0101\cdots 01 \tilde{W_{1}}\tilde{W_{2}}\cdots\tilde{W_{k}},\]
where the alternating sequence $0101\cdots01$ has length $2k$ and
\[\tilde{W_{k}}=W_{k}+1+\sum_{j=1}^{k-1}N_{j}\]
is the sequence obtained by increasing each value in $W_{k}$ by $1+\sum_{j=1}^{k-1}N_{j}$. Then we need to prove the following three facts:
\begin{itemize}
\item The sequence $C$ is uniquely defined by $W_{1},\ldots,W_{k}$,
\item The sequence $C$ avoids $\pi$,
\item The sequence $C$ is an ascent sequence.
\end{itemize}
To show that $C$ is uniquely defined by $W_{1},\ldots,W_{k}$, we break $C$ into $n+2$ subsequences of length $n$ to find the subsequences $\tilde{W_{1}},\tilde{W_{2}},\ldots,\tilde{W_{k}}$. Then $W_{1},\ldots,W_{k}$ can be determined as follows: $W_{1}=\tilde{W_{1}}-1$, then $N_{1}=\max W_{1}$, then $W_{2}=\tilde{W_{2}}-N_{1}-1$, then $N_{1}=\max W_{1}$ and so on until all $W_{i}$ and $N_{i}$ are determined.

Now we will show that $C$ avoids $\pi$. 
We know that each $W_{i}$, and hence each $\tilde{W_{i}}$, avoids $\pi$, as does $01\cdots 01$, so if the pattern $\pi$ appears in $C$, its final element(s) must lie in some $\tilde{W_{i}}$ that does not contain all of its elements. 
But then $\pi$ would decompose as a direct sum, as the element in $\tilde{W_{i}}$ would be greater than all elements not in $\tilde{W_{i}}$. This is a contradiction as we assumed $\pi$ was sum-indecomposable.

Finally we will prove that the sequence $C$ is an ascent sequence. Let $C=n_{1}n_{2}\cdots n_{k^{2}+2k}$ and let $i\in[2,k^{2}+2k]$. Then we need to show that $n_{i}\leq asc(n_{1}\ldots n_{i-1})+1$. If $i\in[2,2k]$ then $n_{i}\in\{0,1\}$, so this is clear. Otherwise, assume $n_{i}$ is an element of $\tilde{W_{h}}$. Then
\begin{align*}asc(n_{1}\ldots n_{i-1})&\geq asc(0101\cdots01\tilde{W_{1}}\cdots\tilde{W_{h-1}})\\
&\geq k+asc(W_{1})+\ldots+asc(W_{h-1})\\
&\geq k+N_{1}+\ldots+N_{h-1}.\end{align*}
Now
\begin{align*}n_{i}\leq\max(W_{h})&=N_{1}+\ldots+N_{h-1}+N_h\\
&\leq N_{1}+\ldots+N_{h-1}+asc(W_h)\\
&\leq N_{1}+\ldots+N_{h-1}+n-1.\end{align*}
Combining these yields $n_{i}\leq asc(n_{1}\ldots n_{i-1})$, so $C$ is an ascent sequence. This completes the proof that $c_{k^2+2k}\geq w_{k}^{k}$.

Finally, we will show that $\mu_{c}=\mu_{w}$. Since $w_{k}\geq c_{k}$, we clearly have $\mu_{w}\geq \mu_{c}$. Moreover,
\[\mu_{c}=\lim_{k\to\infty}c_{k^{2}+2k}^{1/(k^{2}+2k)}\geq \lim_{k\to\infty}w_{k}^{k/(k^{2}+2k)}=\lim_{k\to\infty}w_{k}^{1/(k+2)}=\mu_{w}.\]
\end{proof}

\begin{theorem}
The growth rate $\mu_{c}^{120}$ of $120$-avoiding ascent sequences is equal to the growth rate $\mu_{c}^{201}$ of $201$-avoiding ascents sequences.
\end{theorem}
\begin{proof}
Using Proposition \ref{prop:weak-ascent}, it suffices to show that the growth rates $\mu_{w}^{120}$ and $\mu_{w}^{201}$ of weak ascent sequences are equal. In fact we prove the stronger result that for any $k$, the number $w_{k}^{120}$ of $120$-avoiding weak ascent sequences of length $k$ is equal to the number $w_{k}^{201}$ of $201$-avoiding weak ascent sequences of length $k$. We will show this by a bijection.

Let $n_{1}\cdots n_{k}$ be a 120-avoiding weak ascent sequence. Define the sequence $m_{1}\cdots m_{k}$ by $m_{j}=\max\{n_{1}\cdots n_{k}\}+\min\{n_{1}\cdots n_{k}\}-n_{k+1-j}$. Note that the maximum and minimum values of the sequence do not change under this transformation, so applying this transformation a second time yields the original sequence. Also note that $n_{i}<n_{i+1}$ if and only if $m_{k-i}<m_{k+1-i}$, so the two sequences have the same number of ascents. Hence one is an ascent sequence if and only if the other is an ascent sequence. Finally $n_{i_{1}}$, $n_{i_{2}}$, $n_{i_{3}}$ have the shape $120$ if and only if $m_{k+1-i_{3}}$, $m_{k+1-i_{2}}$, $m_{k+1-i_{1}}$ have the shape $201$. Hence this indeed forms a bijection between $120$-avoiding weak ascent sequences and $201$-avoiding weak ascent sequences. 
\end{proof}

\section{000-avoiding ascent sequences}
This sequence $\{c_n \}$ is given as A202058 in the OEIS \cite{OEIS} to order O$(x^{22})$, and we have extended it to O$(x^{395}).$  We first plotted the ratios of the coefficients $ c_n/c_{n-1}$ against $1/n.$ If one has a pure power-law singularity, such a plot should be linear, with ordinate intercept giving an estimate of the growth constant. The ratio plot (not shown) is clearly diverging as $n \to \infty,$ implying zero radius of convergence. 

\begin{figure}[h!] 
\begin{minipage}[t]{0.45\textwidth} 
\centerline{\includegraphics[width=\textwidth]{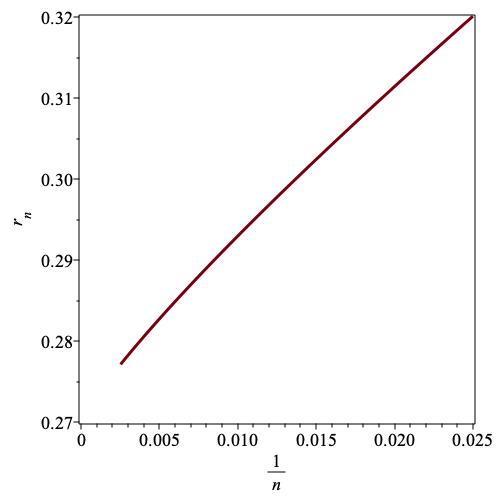}}
%\captionsetup{justification=centering}
\caption{
Plot of ratios of the e.g.f. against $1/n$ for $000$-avoiding ascent sequences.} 
\label{fig:rat1000}
\end{minipage}
\hspace{0.05\textwidth}
\begin{minipage}[t]{0.45\textwidth} 
\centerline {\includegraphics[width=\textwidth]{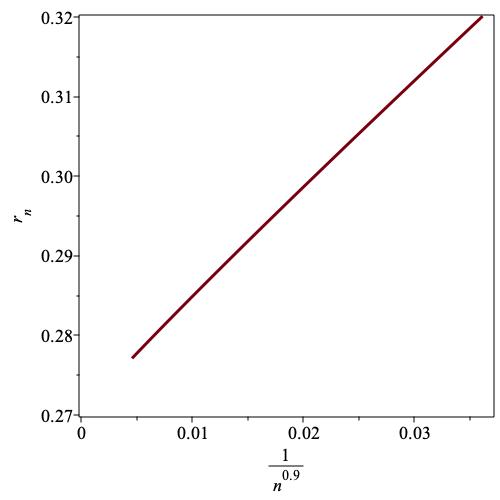}}
% \captionsetup{justification=centering}
\caption{Plot of ratios of the e.g.f. against $1/n^{0.9}$ for  $000$-avoiding ascent sequences. }
\label{fig:rat2000}
\end{minipage}
\end{figure}

This suggest  that one should be looking at the ratios of the exponential generating function (e.g.f.), $r_n=c_n/(n\cdot c_{n-1}), $ which are shown in Fig. \ref{fig:rat1000}. While apparently going to a finite limit as $n \to \infty,$ this displays some curvature. By contrast, if the ratios are plotted against $1/n^{0.9},$ as shown in  Fig. \ref{fig:rat2000}, the plot is virtually linear, and intercepts the ordinate at around 0.271 or 0.272, which is our initial estimate of the growth constant. 

\begin{figure}[h!] 
\begin{minipage}[t]{0.45\textwidth} 
\centerline{
\includegraphics[width=\textwidth]{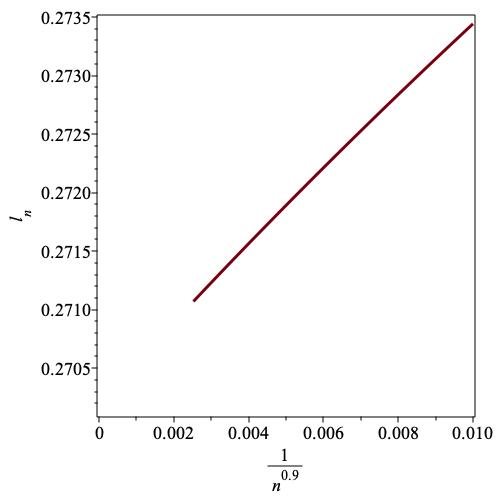}}
%\captionsetup{justification=centering}
\caption{Plot of linear intercepts $l_n$ against $1/n^{0.9}$ for $000$-avoiding ascent sequences.}
 \label{fig:ln000}
\end{minipage}
\hspace{0.05\textwidth}
\begin{minipage}[t]{0.45\textwidth} 
\centerline {\includegraphics[width=\textwidth]{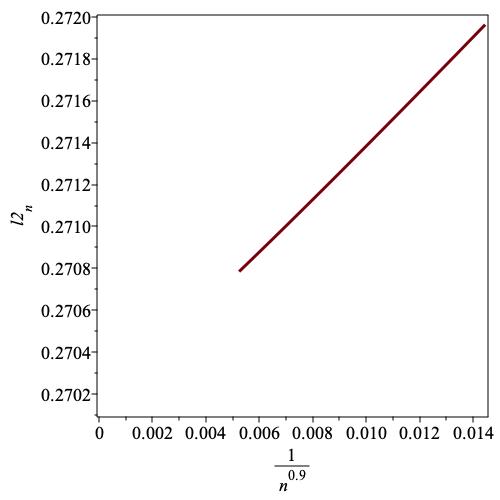}}
% \captionsetup{justification=centering}
\caption{Plot of quadratic intercepts $l2_n$ against $1/n^0.9$ for $000$-avoiding ascent sequences..} 
\label{fig:rat3000}
\end{minipage}
\end{figure}

However linearity against $1/n^{0.9}$ is very close to $1/n$ so perhaps this apparent behaviour is due to the effect of higher-order terms mixed with a O$(1/n)$ term? To eliminate the presumed O$(1/n)$ term, we show in Fig. \ref{fig:ln000} the linear intercepts, $l_n=n\cdot r_n - (n-1) \cdot r_{n-1}$ plotted against $1/n^{0.9}.$ This plot is also linear, and again has ordinate interception between 0.2701 and 0.2702. 

\begin{figure}[h!] 
\begin{minipage}[t]{0.45\textwidth} 
\centerline{
\includegraphics[width=\textwidth]{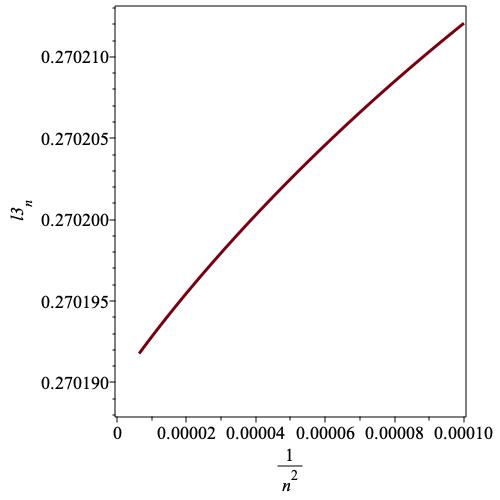}}
%\captionsetup{justification=centering}
\caption{Plot of linear extrapolants of $l2_n$ against $1/n^{2}$ for $000$-avoiding ascent sequences.}
 \label{fig:rat4000}
\end{minipage}
\hspace{0.05\textwidth}
\begin{minipage}[t]{0.45\textwidth} \centerline{\includegraphics[width=\textwidth]{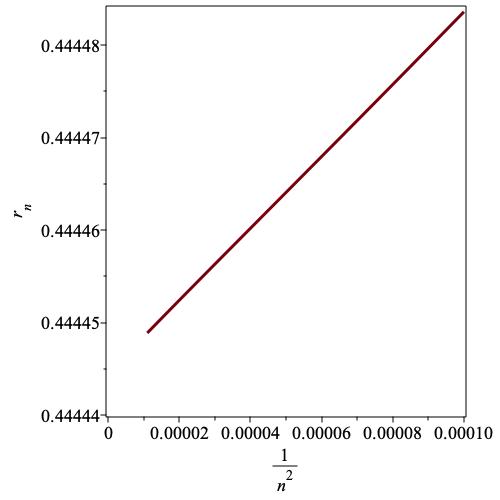}}
%\captionsetup{justification=centering}
\caption{
Plot of linear intercepts of extrapolated Hadamard quotients against $1/n^2$ for $000$-avoiding ascent sequences.} 
\label{fig:hrat}

 \end{minipage}
\end{figure}

We next eliminate terms of O$(1/n^2)$ by constructing the sequence $$l2_n =\frac{n^2\cdot l_n - (n-1)^2 \cdot l_{n-1}}{2n-1} \sim \mu  (1 + {\rm higher\,\,order\,\,terms}).$$ This is shown in Fig. \ref{fig:rat3000}, where $l2_n$ appears linear when plotted against $1/n^{0.9}.$ Since we've eliminated the term O$(1/n)$ and the O$(1/n^2)$
terms in the ratios by constructing the sequence $l2_n,$ the fact that we still have linearity when plotted against $1/n^{0.9}$ implies that there is a term of this form in the ratios. The elimination of the O$(1/n)$ and O$(1/n^2)$ terms in the ratios would not affect such a term, and that is what we are seeing.  The result of extrapolating successive pairs of points by constructing the sequence given by $l3_n=n\cdot l2_n - (n-1) \cdot l2_{n-1}$  is shown in Fig. \ref{fig:rat4000}, from which we estimate $\mu \approx 0.27019 \pm 0.00001.$ The presence of a term of O$(1/n^{0.9})$ implies a stretched-exponential term in the expression for the coefficients, of the form $\mu_1^{n^{\sigma}},$ where $\sigma \approx 0.1.$ We cannot say whether $\sigma =1/10,$ or $1/12,$ but it appears to be in the range $0.06 \le \sigma \le 0.1.$

 Recall that the growth constant for ascent sequences is $6/\pi^2.$ It is interesting to note that $8/(3\pi^2)=0.2701898\ldots.$ Our estimate of 
$\mu$ is in complete agreement with this value, so once again, this is quite suggestive.

Note that all previously solved cases of length-3 avoiding ascent sequences have power-law behaviour. This is the first example where the coefficients grow factorially. Here is a simple argument giving a
 (weak) lower bound to the growth that precludes power-law behaviour. 
 
 Consider an ascent-sequence of length $2n$, the first $n$ terms of which are $0,1,2,\cdots, n-1$. Let the next $n$ terms be any of the $n!$
permutations of $0,1,2,\cdots ,n-1$. Firstly, by construction this sequence  {\em is} an ascent sequence. Secondly, also by construction, it avoids the pattern $000$ as well as the pattern $100.$ Therefore the number of $000$- or $100$-avoiding ascent sequences
of length $2n$ is at least $n!$ So the number of $000$- or $100$-avoiding ascent sequences
of length $n$ is at least $(n/2)!$ 

We can also now compare this sequence with the behaviour of the coefficients of ascent sequences, as they both grow factorially, by constructing the Hadamard quotient of the coefficients of the two sequences. That is to say, the ascent-sequence coefficients are known to grow as $a_n \sim C \cdot n! \cdot (6/\pi^2)^n \cdot \sqrt{n},$ and our best estimate of the behaviour of the coefficients of $000$-avoiding ascent sequences is $c_n \sim D \cdot n! \cdot \mu_1^{n^\sigma}\cdot  \mu^n \cdot n^g.$ The Hadamard quotient  is
$$h_n \equiv \frac{c_n}{a_n} \sim const. \lambda^n \cdot \mu_1^{n^\sigma}  \cdot n^{g-1/2},$$ where $\lambda = \mu \pi^2/6.$ Extrapolation of the ratios $h_n/h_{n-1}$ should give estimates of $\lambda.$ We have eliminated terms of O$(1/n)$ and O$(1/n^2)$ and then plot linear intercepts against $1/n^2$ in Fig. \ref{fig:hrat}. Doing this gives $\lambda \approx 0.444444,$ which implies $\mu = 8/(3\pi^2),$ as suggested earlier, though with significantly greater precision.

The preceding analysis implies that the e.g.f of 000-avoiding ascent sequences behave as $$c_n \sim C \cdot n!\cdot  \mu^n \cdot \mu_1^{n^\sigma} \cdot n^g$$ where we conjecture that $\mu=8/(3\pi^2)$ and $0.06 \le \sigma \le 0.1,$ but we are unable to estimate the value of $\mu_1$ or the exponent $g.$
%To estimate $\sigma$ directly, without assuming the value of $\mu,$ we proceed as we did in our analysis of 120-avoiding ascent sequences in the previous section, and estimate $\sigma$ by using the estimators given in eqns. (\ref{eq:sig1}) and (\ref{eq:sig2}). These are shown in Fig. \ref{fig:sig1000}. It can be seen that both are going to a value around $\sigma-2 \approx -1.9,$ supporting our previous estimate $\sigma \approx 0.1.$

%\centerline{\includegraphics[width=\textwidth,angle=0]{Av25314rat.jpg} }

\section{100-avoiding ascent sequences}
This sequence $\{c_n \}$ is given as A202059 in the OEIS \cite{OEIS} to order O$(x^{23})$, and we have extended this to O$(x^{712}).$  We first plotted the ratios of the coefficients $ c_n/c_{n-1}$ against $1/n.$  As for the 000-avoiding ascent sequences he ratio plot (not shown) is clearly diverging as $n \to \infty,$ implying a zero radius of convergence. We have shown above that $(n/2)!$ is a lower bound for the coefficients of this ascent sequence,
so this result is not surprising.

\begin{figure}[h!] 
\begin{minipage}[t]{0.45\textwidth} 
\centerline{
\includegraphics[width=\textwidth]{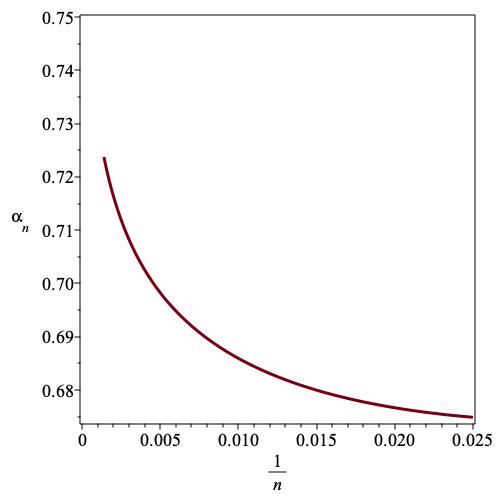}}
%\captionsetup{justification=centering}
\caption{Plot of estimators $\alpha_n$ against $1/n$ for $100$-avoiding ascent sequences.}
 \label{fig:alf100}
\end{minipage}
\hspace{0.05\textwidth}
\begin{minipage}[t]{0.45\textwidth} 
\centerline {\includegraphics[width=\textwidth]{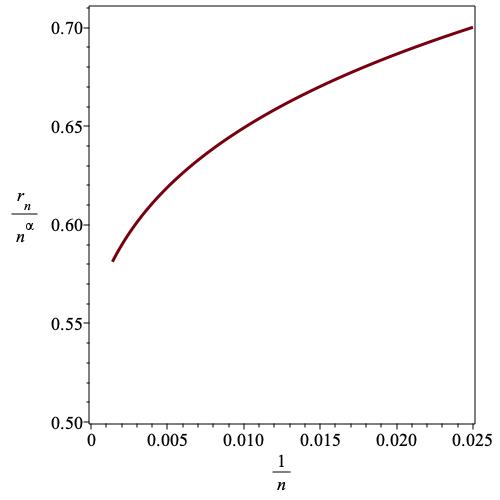}}
% \captionsetup{justification=centering}
\caption{Plot of estimators $\mu\cdot \alpha^{\alpha}$ against $1/n$ for $100$-avoiding ascent sequences.} 
\label{fig:mu100}
\end{minipage}
\end{figure}

Let's assume that the asymptotics are \BE \label{eq:c100}
c_n \sim C \left (\alpha n \right ) ! \mu^n n^g.
\EE
Then
\BE \label{eq:r100}
r_n=\frac{c_n}{c_{n-1}} \sim \alpha^{\alpha}\cdot n^\alpha  \cdot \mu \cdot \left ( 1 + \frac{g}{n} \right ),
\EE
and
\BE \label{eq:s100}
s_n=\frac{r_n}{r_{n-1}} \sim  \left ( 1 + \frac{\alpha}{n} - \frac{g}{n^2}\right ).
\EE
So from eqn. (\ref{eq:s100}) a plot of $s_n$ against $1/n$ should approach ordinate 1 with gradient  $\alpha.$ The term O$(1/n^2)$ can cause some curvature, so we eliminate this by
forming 
\BE \label{eq:tn}
t_n \equiv (n^2s_n-(n-1)^2 s_{n-1})/(2n-1) \sim (1+\alpha/(2n) +o(1/n^2)).
\EE
 We show in Fig. \ref{fig:alf100} estimators of $\alpha$ obtained from the gradient of
the plot of $t_n$ against $1/n.$ We estimate the ordinate to be around $\alpha = 0.75.$

Assuming $\alpha = 3/4,$ we estimate $\mu,$ following eqn. (\ref{eq:r100}), from the ordinate of the plot of $n^{-\alpha}\cdot r_n$ against $1/n.$ Again, we eliminate the O$(1/n^2)$ term as above. The result is shown in
Fig. \ref{fig:mu100}, and we estimate the ordinate to be $0.53 \pm 0.02,$ from which follows $\mu=0.66 \pm 0.02.$ (This analysis assumes $\alpha =3/4.$ A slightly different value of $\alpha$ would make a significant difference
to this estimate.)

An alternative analysis follows by direct fitting. From eqn. (\ref{eq:c100}) we have
 $$\log{c_n} \sim \alpha \cdot n\log{n} +n(\log{\mu} +\alpha\log{\alpha} -\alpha) +(g+1/2)\log{n} + \log(C\sqrt{2\pi\alpha}),$$ where we have used Stirling's approximation for the factorial function.
So fitting $$\log{c_k} = e_1\cdot k\log{k} +e_2 \cdot k +e_3 \cdot \log{k} + e_4,$$ to successive coefficients $\log{c_k}$ with $k=m-2,\,m-1,\, m,\, m+1,$ with $m$ increasing until we run out of known coefficients, we obtain a system of linear equations
from which $e_1,\,e_2,\,e_3,\,e_4$ give estimates of the critical parameters.

\begin{figure}[h!] 
\begin{minipage}[t]{0.45\textwidth} 
\centerline{
\includegraphics[width=\textwidth]{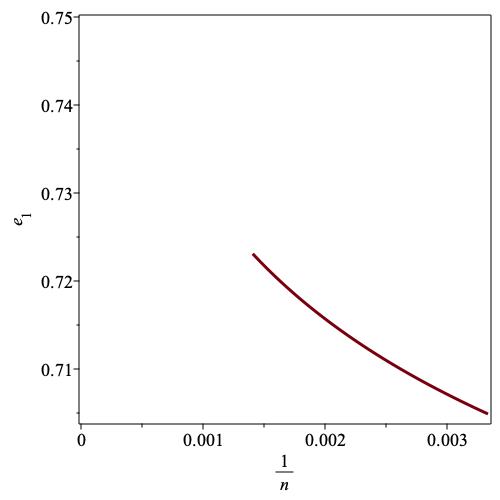}}
%\captionsetup{justification=centering}
\caption{Plot of estimators $e_1=\alpha$ against $1/n$ for $100$-avoiding ascent sequences.}
 \label{fig:e1100}
\end{minipage}
\hspace{0.05\textwidth}
\begin{minipage}[t]{0.45\textwidth} 
\centerline {\includegraphics[width=\textwidth]{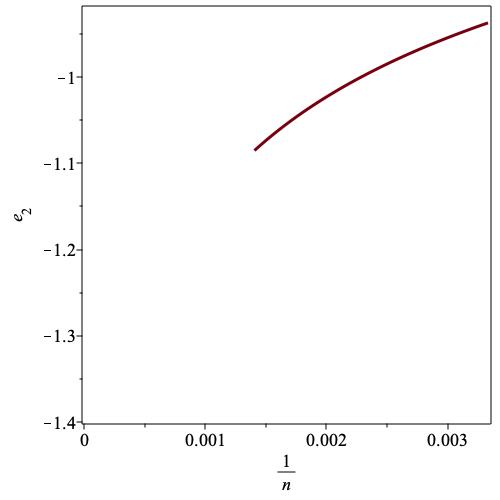}}
% \captionsetup{justification=centering}
\caption{Plot of estimators of $e_2=\log{\mu} +\alpha\log{\alpha} -\alpha$ against $1/n$ for $100$-avoiding ascent sequences.} 
\label{fig:e2100}
\end{minipage}
\end{figure}

In Fig. \ref{fig:e1100} we plot estimates of $e_1=\alpha$ against $1/n.$ Our ratio estimate of $\alpha = 3/4$ is well supported. In Fig. \ref{fig:e2100} we plot estimates of $e_2=\log{\mu} +\alpha\log{\alpha} -\alpha$ against $1/n.$ We estimate
$e_2=-1.35 \pm 0.05,$ from which follows $\mu = 0.68 \pm 0.04,$ in reasonable agreement with the ratio estimate, $\mu =0.66 \pm 0.02.$ We are not confident estimating the other parameters $g$ and $C,$ as they depend sensitively on the values of $\alpha$ and $\mu.$

Accepting that the dominant growth term is $(3n/4)!,$ one can divide the coefficients by this term and analyse the resulting sequence, which hopefully behaves like a conventional power-law singularity. Unfortunately, doing this did
not improve our analysis. The ratio plot of the coefficients still exhibited significant curvature, making it difficult to estimate the various assumed critical parameters. It is possible that the sub-dominant asymptotics is more complicated than we have assumed, but
we have no good idea how to explore this possibility.

We conclude this section with the estimate that the coefficients of 100-avoiding ascent sequences behave as $$ c_n \sim C \left (\alpha n \right ) ! \mu^n n^g,$$ with $\alpha \approx 3/4,$ and $\mu = 0.68 \pm 0.04,$ assuming $\alpha = 3/4$ exactly.
We make no estimate of $C$ or $g.$

\section{110-avoiding ascent sequences}
This sequence $\{c_n \}$ is given as A202060 in the OEIS \cite{OEIS} to order O$(x^{17})$, and we have extended this to O$(x^{42}).$  We have used these exact coefficients to derive 100 further approximate coefficients by the method of series extension \cite{G16}, and briefly described in the Appendix. We first plotted the ratios of the coefficients $ c_n/c_{n-1}$ against $1/n.$ As with the sequence for 100-avoiding ascent sequences, studied in the preceding section, the ratio plot (not shown) is clearly diverging as $n \to \infty,$ implying zero radius of convergence. 

This is not surprising, as a variation of our previous argument gives a lower bound $(n/3)!$ for the growth of the coefficients. Consider an ascent-sequence of length $3n,$ the first $n$ terms of which are $0,1,2,\cdots,n-1$ as is the second block of $n$ terms. Let the next $n$ terms be any of the $n!$ permutations of $n,n+1, n+2, \cdots 2n-1.$ Firstly, by construction this sequence  {\em is} an ascent sequence. Secondly, also by construction, it avoids the pattern $110.$  Therefore the number of $110$-avoiding ascent sequences of length $3n$ is at least $n!$ 
So the number of $110$-avoiding ascent sequences of length $n$ is at least $(n/3)!$.

Our analysis closely parallels that of 100-avoiding ascent sequences, described in the preceding section, and indeed, we find the asymptotic behaviour to be similar, just with a different growth constant.

We first assume that the asymptotics are \BE 
c_n \sim C \left (\alpha n \right ) ! \mu^n n^g.
\EE

\begin{figure}[h!] 
\begin{minipage}[t]{0.45\textwidth} 
\centerline{
\includegraphics[width=\textwidth]{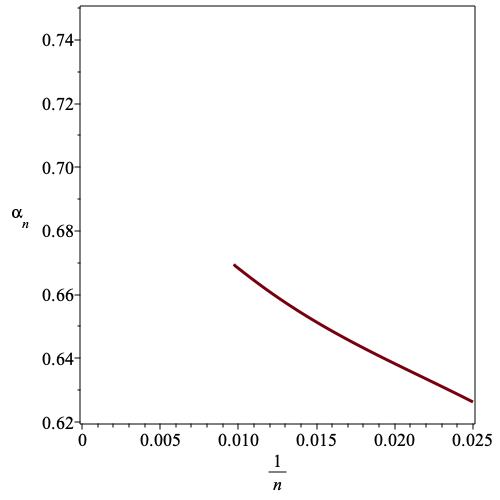}}
%\captionsetup{justification=centering}
\caption{Plot of estimators $\alpha_n$ against $1/n$ for $110$-avoiding ascent sequences.}
 \label{fig:alf110}
\end{minipage}
\hspace{0.05\textwidth}
\begin{minipage}[t]{0.45\textwidth} 
\centerline {\includegraphics[width=\textwidth]{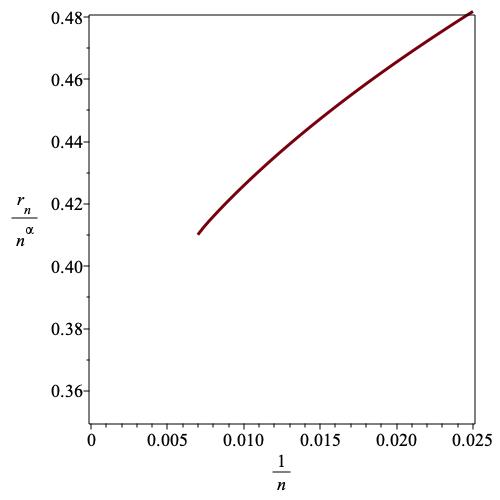}}
% \captionsetup{justification=centering}
\caption{Plot of estimators $\mu\cdot \alpha^{\alpha}$ against $1/n$ for $110$-avoiding ascent sequences.} 
\label{fig:mu110}
\end{minipage}
\end{figure}

So as with the analysis of 100-avoiding ascent sequences,
we  show in Fig. \ref{fig:alf110} estimators of $\alpha$ obtained from the gradient of
the plot of $t_n$ (\ref{eq:tn}) against $1/n.$ We again estimate the ordinate to be around $\alpha = 0.75.$

Assuming $\alpha = 3/4,$ we estimate $\mu,$ following eqn. (\ref{eq:r100}), from the ordinate of the plot of $n^{-\alpha}\cdot r_n$ against $1/n.$ Again, we eliminate the O$(1/n^2)$ term as above. The result is shown in
Fig. \ref{fig:mu110}, and we estimate the ordinate to be $0.36 \pm 0.02,$ from which follows $\mu=0.45 \pm 0.02.$ (This analysis assumes $\alpha =3/4.$ A slightly different value of $\alpha$ would make a significant difference
to this estimate.)

An alternative analysis follows by direct fitting, just as in the previous case.
 
\begin{figure}[h!] 
\begin{minipage}[t]{0.45\textwidth} 
\centerline{
\includegraphics[width=\textwidth]{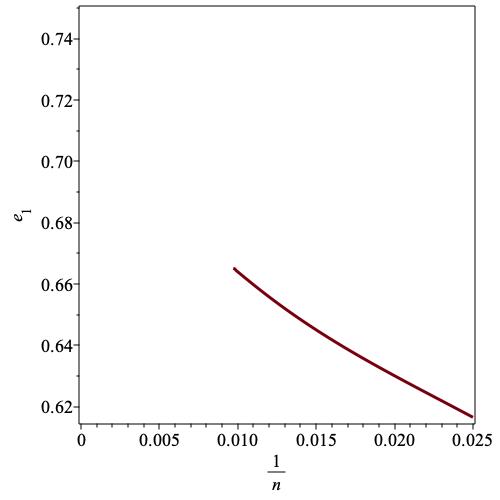}}
%\captionsetup{justification=centering}
\caption{Plot of estimators $e_1=\alpha$ against $1/n$ for $110$-avoiding ascent sequences.}
 \label{fig:e1110}
\end{minipage}
\hspace{0.05\textwidth}
\begin{minipage}[t]{0.45\textwidth} 
\centerline {\includegraphics[width=\textwidth]{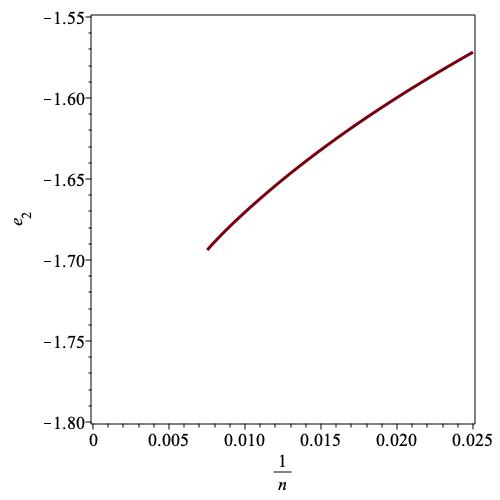}}
% \captionsetup{justification=centering}
\caption{Plot of estimators of $e_2=\log{\mu} +\alpha\log{\alpha} -\alpha$ against $1/n$ for $110$-avoiding ascent sequences.} 
\label{fig:e2110}
\end{minipage}
\end{figure}

In Fig. \ref{fig:e1110} we show estimates of $e_1=\alpha$ plotted against $1/n.$ The ratio estimate of $\alpha = 3/4$ is moderately well supported, though not as well as for the 100-avoiding sequence. However there we had over 700 series coefficients. In Fig. \ref{fig:e2110} we plot estimates of $e_2=\log{\mu} +\alpha\log{\alpha} -\alpha$ against $1/n,$ made by fitting with $\alpha$ assumed to be $3/4.$ We estimate
$e_2=-1.79 \pm 0.05,$ from which follows $\mu = 0.44 \pm 0.02,$ in reasonable agreement with the ratio estimate, $\mu = 0.45 \pm 0.02.$ We are not confident estimating the other parameters $g$ and $C,$ as they depend sensitively on the values of $\alpha$ and $\mu.$

We conclude this section with the estimate that the coefficients of 100-avoiding ascent sequences behave as $$ c_n \sim C \left (\alpha n \right ) ! \mu^n n^g,$$ with $\alpha \approx 3/4,$ and $\mu = 0.44 \pm 0.02,$ assuming $\alpha = 3/4$ exactly.
We make no estimate of $C$ or $g.$

For both this sequence and the 100-avoiding sequence, we conjectured that the dominant growth term is $(3n/4)!.$ We explore this further by considering the Hadamard quotient of the two sequences, which should then have exponential growth.

\begin{figure}[h!] 
\begin{minipage}[t]{0.45\textwidth} 
\centerline{
\includegraphics[width=\textwidth]{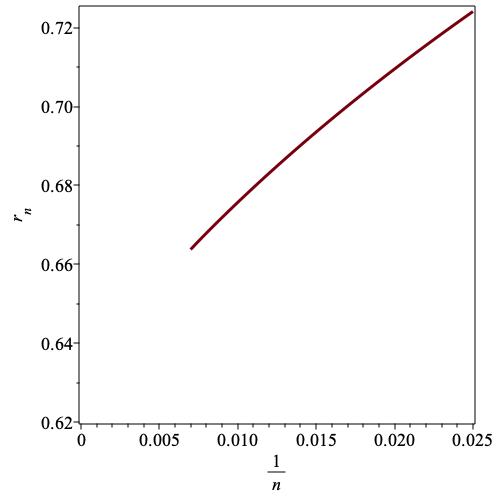}}
%\captionsetup{justification=centering}
\caption{Plot of ratios $r_n$ against $1/n$ for the coefficients $a_n$.}
 \label{fig:c1a}
\end{minipage}
\hspace{0.05\textwidth}
\begin{minipage}[t]{0.45\textwidth} 
\centerline {\includegraphics[width=\textwidth]{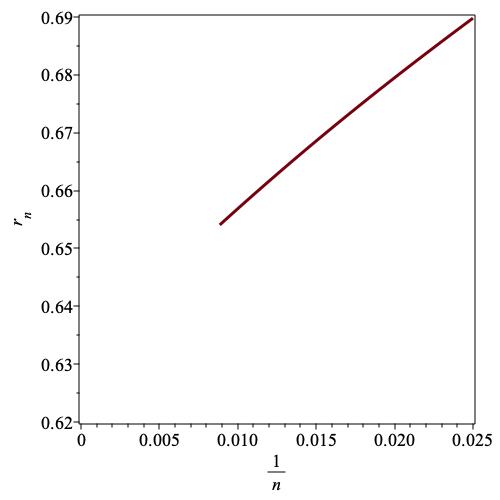}}
% \captionsetup{justification=centering}
\caption{The same plot, with terms O$(1/n^2)$ eliminated.} 
\label{fig:c2a}
\end{minipage}
\end{figure}

Define new coefficients $a_n \equiv c100_n/c110_n,$ which should behave as $D \cdot \mu_{100}/\mu_{110},$ where $D$ is a constant. We study the behaviour of the coefficients $a_n$ by the ratio method, and in Fig \ref{fig:c1a} we show the ratios $r_n=a_n/a_{n-1}$ plotted against $1/n.$ Because of the significant curvature, we eliminate terms of O$(1/n^2),$ which frequently cause this,  and plot the result, also against $1/n,$ in Fig \ref{fig:c2a}. We estimate the limit as $0.625 \pm 0.015.$ From the direct estimates of the growth constants, we find their ratio to be $0.65 \pm 0.07,$ so studying the ratios of the Hadamard quotients gives a more precise estimate of the ratio of the growth constants. The convergence adds support to our belief that the factorial growth is the same for the two sequences.

\section{Conclusion}
We have given a new algorithm to generate length-3 pattern-avoiding ascent sequences, and used this to generate many coefficients for 120-avoiding, 000-avoiding, 100-avoiding and 110-avoiding ascent sequences. 

In each case we have given, conjecturally, the asymptotics of the coefficients. For 120-avoiding ascent sequences we prove that the value of the growth constant is the same as that for 201-avoiding ascent sequences, which is known.
For 000-avoiding ascent sequences we can reasonably confidently conjecture the exact value of the growth constants, and in three cases we have given weak lower bounds that prove that the super-exponential growth conjectured is to be expected. 

These are the only examples
of length-3 pattern-avoiding ascent sequences whose generating function has zero radius of convergence.

\section{Acknowledgements}
AJG would like to thank the ARC Centre of Excellence for Mathematical and Statistical Frontiers (ACEMS) for support. This research was undertaken using the Research Computing Services facilities hosted at the University of Melbourne.

%\end{document}
\section{Appendix}
\section{Series analysis}
The method of series analysis  has, for many years, been a powerful tool in the study of a variety of problems in statistical mechanics, combinatorics, fluid mechanics and computer science. In essence, the problem is the following: Given the first $N$ coefficients of the series expansion of some function, (where $N$ is typically as low as 5 or 6, or as high as 100,000 or more), determine the asymptotic form of the coefficients, subject to some underlying assumption about the asymptotic form, or, equivalently, the nature of the singularity of the function.

\subsection{Ratio Method}
The ratio method was perhaps the earliest systematic
method of series analysis employed,
and is still the most useful method when only a small number of terms are known.
If we have a power-law singularity, so that $f(z) =\sum c_n z^n \sim C(1- z/z_c)^{-\gamma},$
it follows that the {\it ratio} of successive terms
\begin{equation} \label{ratios}
r_n = \frac{c_n}{c_{n-1}}=\frac{1}{z_c}\left (1 + \frac{\gamma -1}{n} + {\rm o}(\frac{1}{n})\right ).
\end{equation}
 It is then natural to plot the successive ratios $r_n$ against $1/n.$
If the correction terms ${\rm o}(\frac{1}{n})$ can be ignored\footnote{For a purely algebraic singularity with no confluent terms, the correction term will be ${\rm O}(\frac{1}{n^2}).$}, such a plot will be linear,
with gradient $\frac{\gamma-1}{z_c},$ and intercept $\mu=1/z_c$ at $1/n = 0.$

\subsection{Functions with non-power-law singularities.}\label{sec:nonalg}

A number of solved, and unsolved problems that arise in lattice critical phenomena and algebraic combinatorics have coefficients with a more complex asymptotic form, with a sub-dominant  term ${\rm O}(\mu_1^{n^\sigma})$ as well as a power-law term ${\rm O}(n^g).$  Perhaps the best-known example of this sort of behaviour is the number of partitions of the integers -- though in that case the leading exponential growth term $\mu^n$ is absent (or equivalently $\mu=1$).
The  form of the coefficients $c_n$ in the general case is
\BE \label{stretch}
c_n \sim C \cdot \mu^n \cdot \mu_1^{n^\sigma} \cdot n^g.
\EE

An example from combinatorics is given by Dyck paths enumerated not just by
length, but also by height (defined to be the maximum vertical distance of the
path from the horizontal axis). Let $d_{n,h}$ be the number of Dyck paths of
length $2n$ and height $h.$ The OGF is then\footnote{One of us (AJG) posed this problem at an Oberwolfach meeting in March 2014. Within 24 hours Brendan McKay produced this solution.}
\BE
D(x,y) = \sum_{n,h} d_{n,h}x^{2n}y^h,\,\,{\rm and} \,\, [x^{2n}]D(x,y)=\sum_{h=1}^n d_{n,h}y^h.
\EE
For $y < 1$
let $A=2^{5/3}\pi^{5/6}/\sqrt{3},$ $E=3\left ( \frac{\pi}{2} \right )^{2/3}$ and $r=-\log{y}.$
Then one finds that $[x^{2n}]D(x,y)$ is given by eqn. (\ref{stretch}) with $C=\frac{1-y}{y^2}r^{1/3}A,$ $\mu=4,$ $\mu_1=\exp(-Er^{2/3}),$ $ \sigma=1/3,$ and $g=-5/6.$

Applying the ratio method to such singularities requires some significant changes. These were first developed in \cite{G15}, where further details and more examples can be found. 
In the next subsection we give a summary, including as much detail as is needed for our analysis.

\subsection{Ratio method for stretched-exponential singularities.} \label{sec:se}

If 

\BE \label{eq:an}
b_n \sim C \cdot \mu^n \cdot \mu_1^{n^\sigma} \cdot n^g,
\EE
 then the ratio of successive coefficients $r_n = b_n/b_{n-1},$ is
\begin{multline} \label{eq:rna}
r_n = \mu \left (1 + \frac{\sigma \log \mu_1}{n^{1-\sigma}} + \frac{g}{n} + \frac{\sigma^2 \log^2 \mu_1}{2n^{2-2\sigma}} + \frac {(\sigma-\sigma^2)\log \mu_1+2g\sigma \log \mu_1}{2n^{2-\sigma}} \right . \\
 \left . {}+ \frac{\sigma^3 \log^3 \mu_1}{6n^{3-3\sigma}} +{\rm O}(n^{2\sigma-3}) + {\rm O}(n^{-2}) \right ).
%\notag
\end{multline}

It is usually the case that $\sigma$ takes the simple values $1/2,$ $1/3,$ $1/4$ etc.\footnote{In statistical mechanical models, the value of the exponent $\sigma$ is simply related to the fractal dimension $d_f$  of the object through $\sigma = 1/(1+d_f).$ }. If these asymptotics arise as the irregular singular point of a D-finite ODE, the exponent must be of the form $1/\rho,$ where $\rho$ is a positive integer.

The presence of the term O$(\frac{1}{n^{1-\sigma}})$ in the expression for the ratios above means that a ratio plot against $1/n$ will display curvature, which can be usually be removed by 
 plotting the ratios against $1/n^{1-\sigma}.$ 
   
 Unfortunately the observation that a ratio plot against $1/n^{1-\sigma}$ will linearise the plot does not provide a sufficiently precise method to estimate the value of $\sigma.$  One can usually distinguish between, say, $\sigma=1/2$ and $\sigma = 1/3$ in this way, but one cannot be much more precise than that. However, one can extend the ratio method to provide direct estimates for the value of $\sigma.$

From (\ref{eq:rna}), one sees that 
\BE \label{eq:rsigma}
(r_n/\mu-1) = \sigma \log{\mu_1} \cdot n^{\sigma-1} + O\left ( \frac{1}{n} \right ).
\EE
 Accordingly, a plot of $\log(r_n/\mu-1)$ versus $\log{n}$ should be linear, with gradient $\sigma-1.$ We would expect an estimate of $\sigma$ close to that which linearised the ratio plot.

This log-log plot will usually be visually linear, but the local gradients are changing slowly as $n$ increases. It is therefore worthwhile extrapolating the local gradients. To do this, from (\ref{eq:rsigma}), we form the estimators
\BE \label{eqn:sig1}
{\tilde \sigma}_n = 1+ \frac{\log |r_n/\mu-1|-\log |r_{n-1}/\mu-1|}{\log{n}-\log(n-1)}.
\EE
This can be extrapolated against $1/n^\sigma,$ using any approximate value of $\sigma.$

A second estimator of $\sigma$ follows from eqn. (\ref{eq:an}). Define $$c_n \equiv \log(b_n/\mu^n) \sim \log{C} + \log{\mu_1}\cdot n^\sigma + g\cdot \log{n},$$  then setting
\BE \label{eq:sigma2}
d_n \equiv c_n-c_{n-1} \sim \sigma \log{\mu_1}\cdot n^{\sigma-1} + g/n,
\EE
 a log-log plot of $d_n$ against $n$ should be linear with gradient $\sigma-1.$ Note that if $\sigma$ is closer to zero than to 1, there is likely to be some competition between the two terms in the  expansion.

This way of estimating $\sigma$ requires knowledge of, or at worst a very precise estimate of, the growth constant $\mu.$ While $\mu$ is exactly known in some cases, more generally $\mu$ is not known, and must be estimated, along with all the other critical parameters. In order to estimate $\sigma$ without knowing $\mu,$ we can use one (or both) of the following estimators:

From eqn. (\ref{eq:rna}), it follows that 
\BE \label{eq:sig1}
r_{\sigma_n} \equiv \frac{r_n}{r_{n-1}} \sim 1 + \frac{(\sigma-1)\log{\mu_1}}{n^{2-\sigma}} + {\rm O}(1/n^2),
\EE
so $\sigma$ can be estimated from a plot of $\log(r_{\sigma_n}-1)$ against $\log{n},$ which should have gradient $\sigma-2.$ Again, the local gradients can be calculated and plotted against $1/n^\sigma,$ using any approximate value of $\sigma.$

Another estimator of $\sigma$ when $\mu$ is not known follows from eqn. (\ref{eq:an}),
\BE \label{eq:sig2}
a_{\sigma_n} \equiv \frac{b_n^{1/n}}{b_{n-1}^{1/(n-1)}} \sim 1 + \frac{(\sigma-1)\log{\mu_1}}{n^{2-\sigma}} + {\rm O}(1/n^2),
\EE
so again $\sigma$ can be estimated from a plot of $\log(a_{\sigma_n}-1)$ against $\log{n}.$ Again, estimates of $\sigma$ are found by extrapolating the local gradient  against $1/n^\sigma.$

While these two estimators are equal to leading order, they differ in their higher-order terms. Which of the two is more informative seems to vary from problem to problem. However, we generally use both.

From eqn. (\ref{eq:rna}), if we know (or conjecture) $\mu$ and $\sigma,$ we can use this to estimate $\mu_1,$ as
\BE \label{eqn:mu1}
\left ( \frac{r_n}{\mu} - 1\right )\cdot n^{1-\sigma} \sim \sigma \cdot \log(\mu_1).
\EE

\section{Differential approximants}
\label{ana:da}

The generating
functions  of some problems in enumerative combinatorics are sometimes algebraic, sometimes D-finite, 
sometimes differentially algebraic, and sometimes transcendentally transcendental.
The not infrequent occurrence of D-finite solutions was the origin of the method of {\em differential approximants}, a very successful method of series analysis for power-law singularities \cite{G89}.

The basic idea is to approximate a generating function $F(z)$ by solutions
of differential equations with polynomial coefficients. That is to say, by D-finite ODEs. The singular behaviour
of such ODEs is  well documented
(see e.g. \cite{Ince27}), and the singular points and
exponents are readily calculated from the ODE. 

The key point for series analysis is that even if {\em globally} the function is not describable by a solution
of such a linear ODE (as is frequently the case) one expects that
{\em locally,} in the
vicinity of the (physical) critical points, the generating
function is still well-approximated by a solution of a linear ODE, when the singularity is a generic power law.

An $M^{th}$-order differential approximant (DA) to a function $F(z)$  is formed by matching
the coefficients in the polynomials $Q_k(z)$ and $P(z)$ of degree $N_k$ and $L$, respectively,
so that the formal solution of the $M^{th}$-order inhomogeneous ordinary differential equation
\BE \label{eq:ana_DA}
\sum_{k=0}^M Q_{k}(z)(z\frac{{\rm d}}{{\rm d}z})^k \tilde{F}(z) = P(z)
\EE
agrees with the first $N=L+\sum_k (N_k+1)$ series coefficients of $F(z)$. 

Constructing such ODEs only involves
solving systems of linear equations. The function
$\tilde{F}(z)$ thus agrees with the power series expansion of the (generally unknown)
function $F(z)$ up to the first $N$ series expansion coefficients.

From the theory of ODEs \cite{Ince27}, the singularities of $\tilde{F}(z)$ are approximated by zeros
$z_i, \,\, i=1, \ldots , N_M$ of $Q_M(z),$ and the
associated critical exponents $\gamma_i$ are estimated from the indicial equation. If there is only a single root at $z_i$  this is just
\BE \label{eq:ana_indeq1}
\gamma_i=M-1-\frac{Q_{M-1}(z_i)}{z_iQ_M ' (z_i)}.
\EE
Estimates of the critical amplitude $C$ are rather more difficult to make, involving the integration of the differential approximant. For that reason the simple ratio method approach to estimating critical amplitudes is often used, whenever possible taking into account higher-order asymptotic terms \cite{GJ09}.

Details as to which approximants should be used and how the estimates from many approximants are averaged to give a single estimate are given in \cite{GJ09}. Examples of the application of the method can be found in \cite{G15}.

 In this work, none of the four series that we analyse are appropriate for analysis by the method of differential approximants, however we describe the method as it underlies the idea of series extension, as described in the next section.
\section{Coefficient prediction}
\label{pred}
In \cite{G16} we showed that the ratio method and the method of differential approximants  work serendipitously together in many cases, even when one has stretched exponential behaviour, in which case neither method works particularly well in unmodified form. 

To be more precise, the method of differential approximants (DAs)  produces ODEs which, by construction, have solutions whose series expansions agree term by term with the known coefficients used in their construction. Clearly, such ODEs implicitly define {\em all}  coefficients in the generating function, but if $N$ terms are used in the construction of the ODE, all terms of order $z^{N}$ and beyond will be approximate, unless the exact ODE is discovered, in which case the problem is solved, without recourse to approximate methods.

What we have found is that it is useful to construct a number of DAs that use all available coefficients, and then use these to predict subsequent coefficients. Not surprisingly, if this is done for a large number of approximants, it is found that the predicted coefficients of the term of order $z^n,$ where $n > N,$ agree for the first $k(n)$ digits, where $k$ is a decreasing function of $n.$ We take as the predicted coefficients the mean of those produced by the various DAs, with outliers excluded, and as a measure of accuracy we take the number of digits for which the predicted coefficients agree, or the standard deviation. These two measures of uncertainty are usually in good agreement.

Now it makes no logical sense to use the approximate coefficients as input to the method of differential approximants, as we have used the DAs to obtain these coefficients. However there is no logical objection to using the ({\em approximate}) predicted coefficients as input to the ratio method. Indeed, as the ratio method, in its most primitive form, looks at a graphical plot of the ratios, an accuracy of 1 part in $10^4$ or $10^5$ is sufficient, as errors of this magnitude are graphically unobservable. 

The DAs use all the information in the coefficients, and are sensitive to even quite small errors in the coefficients. As an example, in a recent study of some self-avoiding walk series, an error was detected in the twentieth significant digit in a new coefficient, as the DAs were much better converged without the last, new, coefficient. The DAs also require high numerical precision in their calculation. In favourable circumstances, they can give remarkably precise estimates of critical points and critical exponents, by which we mean up to or even beyond 20 significant digits in some cases. Surprisingly perhaps, this can be the case even when the underlying ODE is not D-finite. Of course, the singularity must be of the assumed power-law form.

Ratio methods, and direct fitting methods, by contrast are much more robust. The sort of small error that affects the convergence of DAs would not affect the behaviour of the ratios, or their extrapolants, and would thus be invisible to them. As a consequence, approximate coefficients are just as good as the correct coefficients in such applications, provided they are accurate enough. We re-emphasise that, in the generic situation, ratio type methods will rarely give the level of precision in estimating critical parameters that DAs can give. By contrast, the behaviour of ratios can more clearly reveal features of the asymptotics, such as the fact that a singularity is not of power-law type. This is revealed, for example, by curvature of the ratio plots \cite{G15}.

In practice we find, not surprisingly, that the more exact terms we know, the greater is the number of predicted terms, or ratios that can be predicted. 

In this study, we have extended the sequence of coefficients for the generating functions of the two shorter series, 110-avoiding and 120-avoiding ascent sequences by 100 and 200 terms respectively, and have analysed the resulting series by ratio methods.

\end{document}